\documentclass[preprint,12pt]{elsarticle}

\usepackage{graphicx}
\usepackage{color}
\usepackage{subfigure}
\usepackage{amssymb}
\usepackage{amsthm}
\usepackage{amsfonts}

\newtheorem{thm}{Theorem}
\newtheorem{lem}[thm]{Lemma}
\newtheorem{cor}[thm]{Corollary}
\newtheorem{prop}[thm]{Proposition}
\newdefinition{rmk}{Remark}
\newdefinition{exam}{Example}
\newproof{pf}{Proof}
\newproof{pot}{Proof of Theorem \ref{theorem5.1.}}
\newproof{pot1}{Proof of Theorem \ref{theorem4.1.}}

\journal{***}

\begin{document}

\begin{frontmatter}


\title{Threshold dynamics and ergodicity of an SIRS epidemic model
with Markovian switching}

\author[dl]{Dan~Li}
\ead{lidanwhy@163.com}
\author[dl]{Shengqiang~Liu\corref{cor1}}
\ead{sqliu@hit.edu.cn}
\author[jac]{Jing'an~Cui}
\ead{cuijingan@bucea.edu.cn}

\cortext[cor1]{Corresponding author}

\address[dl]{Department of Mathematics, Harbin Institute of Technology, Harbin, Heilongjiang, 150001,
P.R. China}
\address[jac]{School of Science, Beijing University of Civil
  Engineering and Architecture, Beijing, 100044, P.R. China}

\begin{abstract}
This paper studies the spread dynamics
of a stochastic SIRS epidemic model with nonlinear incidence and varying population
size, which is formulated as
a piecewise deterministic Markov process.
A threshold dynamic determined by the basic reproduction
number $\mathcal{R}_{0}$ is established:
the disease can be eradicated almost surely if $\mathcal{R}_{0}<1$,
while the disease persists almost surely if $\mathcal{R}_{0}>1$.
The existing method
for analyzing ergodic behavior of population systems
has been generalized. The modified method weakens the required conditions
and has no limitations
for both the number of environmental
regimes and the dimension of the considered system.
When $\mathcal{R}_{0}>1$,
the existence of a stationary probability measure
is obtained. Furthermore, with the modified method,
the global attractivity of the $\Omega$-limit set of the system
and the convergence in total variation to the stationary measure are
both demonstrated under a mild extra condition.

\end{abstract}


\end{frontmatter}

\section{Introduction}
\label{}

Since the seminal work of Kermack and McKendrick \cite{Kermack-McKendrick1927},
mathematical models have become important tools for
understanding the spread and control of infectious diseases.
In the real ecological systems,
the population dynamics are usually
influenced by a random switching in the external environments.
For example, the disease transmission rate $\beta$ in epidemic models was modified
for meteorological factors because survivals
and infectivity of many viruses and bacteria are better in damp conditions with little
ultraviolet light \cite{Arundel1986,Minhaz2010,Keeling2008};
The population growth rates
and the environmental capacities usually
fluctuate with the change of food resource abundance,
which is in turn dependent upon unpredictable
rainfall fluctuations largely \cite{Dexter2003}.
In the literatures,
it is a popular way that
the random switching of environmental regimes
is characterized by the continuous-time Markov chain
with values in a finite state space, which drives
the changes of the main parameters of
population models
with state switchings of the Markov chain \cite{Liumeng2010,Gray2012,Bao2016Shao,Jiang2015,Jiang2016,Yingang2011,Yingang2014}.
The resulting dynamical systems then become
regime-switching differential equations, i.e.,
the piecewise-deterministic Markov process.
Therefore, based on the biological system subject to stochastic
environmental conditions,
epidemic models with deterministic
parameters are unlikely to be realistic,
and it is significant to investigate the
effect of the random switching of environmental regimes
on the spread dynamics of the disease
in the host population.

However, despite the potential importance of the environmental noise,
it has received relatively little attention in the epidemiology literatures.
Gray, Greenhalgh and Mao et al. \cite{Gray2012}
are the first ones to propose a piecewise deterministic SIS
epidemic model with Markovian switching.
By the modeling techniques of time discretization or branching processes,
Bac\"{a}er et al. \cite{Bacar2013,Bacar2014} obtained
very well results concerning the
definition of the suitable basic
reproduction number for the epidemic model
with random switching of environmental regimes.
Recently, Hieu and Du et al. \cite{HIEU2015} and
Greenhalgh, Liang and Mao \cite{Gray2016} studied
the effect of Markovian switching on the deterministic SIRS epidemic
models respectively, both of which are special cases
of the classic epidemic model as follows
\begin{eqnarray}\label{0.2}
\left\{
  \begin{array}{l}
 d S(t)/dt= \Lambda-\mu S(t)+\lambda R(t)-\beta S(t)I(t),   \\ [+6pt]
 d I(t)/dt= \beta S(t)I(t)-\big(\mu +\alpha +\delta \big)I(t), \\ [+6pt]
 d R(t)/dt= \delta I(t)-\big(\mu +\lambda \big)R(t),
\end{array} \right.
\end{eqnarray}
which was presented by Anderson and May \cite{Anderson1979}
to study the influence of infectious diseases
(caused by viruses or bacteria, etc.) on the density of host populations.
In this model, $S(t)$, $I(t)$ and $R(t)$ are respectively the number (or density) of
susceptible, infectious and recovered individuals at time $t$, and
all the parameters are positive.
$\Lambda$ represents a constant recruitment of new susceptibles;
$\mu$ is the per-capita natural mortality rate;
$\beta$ is the transmission rate (i.e., effective per capita contact rate of infective individuals)
equal to the product of the contact rate and transmission probability;
$\alpha$ is the mortality caused by the disease;
$\delta$ is the per-capita recovery rate of infected individuals;
the recovered hosts are initially immune, but this immunity
can be lost at a rate $\lambda$.

In the existing literatures,
there are few researches on the
ergodicity of the stochastic epidemic model,
which is formulated as the piecewise deterministic Markov process.
To the best of our knowledge,
Hieu and Du et al. \cite{HIEU2015} first studied the
ergodicity of the regime-switching epidemic model,
which was the stochastic edition
of a special case of system (\ref{0.2}).
They obtained the threshold between the extinction
and persistence of the disease.
Using the method of references \cite{Yingang2011,Yingang2014}, they described
completely the $\Omega$-limit set of all positive solutions of the model
and established the sufficient condition ensuring that
the instantaneous measure converges to the stationary measure in total variation.
However, the mortality caused by the disease was not taken into account,
which led to the variable population size
and then resulted in the challenge of reducing the
dimension of the considered system.
Moreover, they assumed that there were only two environmental regimes.

In addition, the incidence function
plays an important role
for ensuring that the epidemic model does give a reasonable
description of the disease dynamics \cite{Capa1993,Levin1989}.
It is traditionally assumed that the incidence rate
of disease spread is bilinear
with respect to the number
of susceptible individuals $S(t)$ and
the number of infective individuals $I(t)$, e.g., $\beta S I$
as in the system (\ref{0.2}).
Actually, it is generally difficult to get the details
of transmission of infectious diseases because they may vary
with the ambient conditions.
Moreover, with the general incidence function, the data themselves
may flexibly decide the function form of incidence rates in practice \cite{Xia-Grenfell2005}.
Therefore, in this paper, we will choose the general nonlinear incidence rate of
the form $\beta S G(I)$ to enable the model to be more realistic
and have wider application.

Motivated by the facts mentioned above, in this paper, we will study
the following stochastic SIRS epidemic model
with the general (including both linear and nonlinear) incidence rate:
\begin{eqnarray}\label{2.2}
\left\{
  \begin{array}{l}
d S(t)/dt=\Lambda-\mu S(t)+\lambda R(t)-\beta_{r(t)} S(t) G(I(t)),   \\ [+6pt]
d I(t)/dt=\beta_{r(t)} S(t) G(I(t))-\big(\mu + \alpha +\delta \big)I(t), \\ [+6pt]
d R(t)/dt=\delta I(t)-\big(\mu +\lambda\big)R(t),
\end{array} \right.
\end{eqnarray}
where $G(\cdot)$ is a general function and
the transmission rate $\beta$ is drived
by a homogeneous continuous-time Markov chain
$\{ r(t), t \geq 0\}$ taking values in
a finite state space $\mathcal{M}=\{1,2,\ldots, E\}$
representing different environments.
As in \cite{Mao1999}, we let the Markov chain $r(t)$ be generated by the
transition rate matrix $Q=(q_{e,e'})_{E\times E}$, i.e.,
\begin{eqnarray*}
\mathbb{P}\{r(t+\Delta t)=e' | r(t)=e\}=
\left\{
  \begin{array}{l}
 q_{e,e'}\Delta t +o(\Delta t),\ \ \ \ \ \ \ \ \  \mbox{ if } \ e \neq e',   \\ [+6pt]
 1+q_{e,e'}\Delta t +o(\Delta t),\ \ \ \ \mbox{ if } \ e = e',
\end{array} \right.
\end{eqnarray*}
where $\Delta t>0$ represents a small time increment.
Here $q_{e,e'}$ is the transition rate from state $e$ to state $e'$,
and $q_{e,e'}\geq 0$ if $e\neq e'$ while
$q_{e,e}=-\sum_{e'\neq e}q_{e,e'}$.

The main aim of this paper is to investigate the extinction
and persistence of the disease of system (\ref{2.2})
and determine the threshold between them.
In the case of the disease persistence, we will
analyze the ergodic behavior of system (\ref{2.2}).
By skilled techniques,
authors in \cite{Yingang2011,Yingang2014} applied the stochastic stability theory
of Markov processes in \cite{Meyn-I,Meyn-II,Meyn-III,Meyn-book}
to two-dimensional Kolmogorov systems of competitive type switching between two
environmental states.
They obtained the sufficient conditions for
the global attractivity of the $\Omega$-limit set of the system
and the convergence of the instantaneous
measure to the stationary measure in total variation.
However, the method used in \cite{Yingang2011,Yingang2014}
is not applicable for our model (\ref{2.2}).
In this paper, we will generalize the method
to analyze ergodic behavior of the ecological systems
from two environmental regimes to any finite ones.
By the theory in \cite{Jurdjevic1997}, we will weaken the condition to ensure the
global attractivity of the $\Omega$-limit set and the convergence of
instantaneous distribution to the stationary distribution.
Moreover, the modified techniques avoid the limitations
for both the number of environmental
regimes and the dimension of the considered system.
For instance, the method used by \cite{Yingang2011,Yingang2014}
does not deal with the case where the dimension of the system
is higher than the number of environmental states of Markov chain,
while the modified method can do.

It is worth mentioning that the modified method can be applied to
the case where
all parameters of the model considered in \cite{Gray2012,HIEU2015,Gray2016}
were derived by the stochastic process $r(t)$,
however, here we only allow the transmission rate $\beta$ of model (\ref{2.2}) to be
disturbed because in reality it may be more sensitive to environmental
fluctuations than other parameters of model (\ref{0.2})
for human populations. In particular,
by the new method, under weaker conditions we can directly extend the results
in \cite{Yingang2011,Yingang2014,HIEU2015}
from two environmental states to any finite ones.

This paper is organized as follows.
Section \ref{Preliminaries-and-main results} introduces some preliminaries used in the later parts and illustrates our main results.
In this section, we investigate the extinction and the persistence in the time mean of the
disease and establish the threshold between them. Under the mild conditions,
we describe completely the $\Omega$-limit set of all positive
solutions to the model (\ref{2.2}),
and prove the global attractivity of the $\Omega$-limit set
and the stochastic stability of the unique invariant measure.
In Section \ref{Proofs of main results}, we give
the proofs of the main results in details.
Finally we conclude this paper with further remarks in Section \ref{Conclusion}.

\section{Preliminaries and main results}
\label{Preliminaries-and-main results}

\subsection{Preliminaries}

In this paper, unless otherwise specified,
let $(\widetilde{\Omega}, \mathcal{F}, \mathbb{P})$ be
a complete probability space with a filtration $\{\mathcal{F}_{t}\}_{t\geq 0}$
satisfying the usual conditions (i.e., it is right continuous and
increasing while $\mathcal{F}_{0}$ contains all $\mathbb{P}$-null sets).
Note that we here use $\widetilde{\Omega}$ instead of the usual $\Omega$
to denote the sample space, still denote by $\omega$ an element of $\widetilde{\Omega}$,
and reserve the notation $\Omega$ for the notion of omega-limit set to avoid notional conflict.
For any initial value $z_{0}=(S(0),I(0),R(0))\in \mathbb{R}^{3}_{+}$
with $\mathbb{R}^{3}_{+}=\big\{x\in \mathbb{R}^{3}: x_{i}>0, i=1, 2, 3 \big\}$,
we denote by $z(t,\omega,z_{0})$ $=$ $(S(t,\omega,z_{0})$, $I(t,\omega,z_{0})$,
$R(t,\omega,z_{0}))$
the solution to system (\ref{2.2}) at time $t$, starting in $z_{0}$
(or $z(t,z_{0})=(S(t,z_{0}),I(t,z_{0}),R(t,z_{0}))$, $z(t)=(S(t),I(t),R(t))$
whenever there is no ambiguity).
If $x\in \mathbb{R}^{3}$ and $\epsilon>0$, the open ball $B(x,\epsilon)$
with center at $x$ and radius $\epsilon$ is defined to be the set of
all $y\in \mathbb{R}^{3}$ such that $|y-x|< \epsilon$.
Denote $\mathbb{N}=\{1, 2, \ldots\}$ and
$\mathbb{N}_{0}=\{0,1, 2, \ldots\}$.
For any constant
sequence $\{c(e):e \in \mathcal{M}\}$, define
$c^{M}=\max_{e\in \mathcal{M}}\{c(e)\}$.

Assume that the Markov chain $r(t)$ is irreducible.
Under this condition, the
Markov chain has a unique stationary probability distribution
$\pi=(\pi_{1}, \ldots, \pi_{E})\in \mathbb{R}^{1\times E}$,
which can be determined by solving
the following linear equation
$\pi Q=0$ subject to
$\sum_{e=1}^{E}\pi_{e}=1$, and $\pi_{e}>0, \ \forall e \in \mathcal{M}$.
Let
$$0=\tau_{0}<\tau_{1}<\tau_{2}<\cdots < \tau_{n}<\cdots$$
be jump times of the Markov chain $r(t)$, and denote by
$$\sigma_{1}=\tau_{1}-\tau_{0}, \ \ \sigma_{2}=\tau_{2}-\tau_{1}, \ \ \ldots,
\ \ \sigma_{n}=\tau_{n}-\tau_{n-1}, \ \ \ldots.$$
holding time between adjacent two jumps. Furthermore,
we set
$$\mathcal{F}_{0}^{n}=\sigma(\tau_{k}, k\leq n)\ \ \mathrm{and} \ \
\mathcal{F}_{n}^{\infty}=\sigma(\tau_{k}-\tau_{n}, k> n)$$
for all $n=0,1,\cdots $, then for each $n$, $\mathcal{F}_{0}^{n}$
is independent of $\mathcal{F}_{n}^{\infty}$.

Throughout this paper, we further assume that
\begin{enumerate}
\item [(\textbf{H1})] $G(\cdot):\mathbb{R}_{+}\rightarrow \mathbb{R}_{+}$
is a twice-differentiable function, $G(0)=0$ and $0< G(I)\leq IG'(0)$ holds for all $I>0$;

\item [(\textbf{H2})] $g(x)$ is Lipschitz on $[0, \Lambda/ \mu]$, namely,
there exists a constant $\vartheta >0$, such that
$$|g(x_{1})-g(x_{2})|\leq \vartheta |x_{1}-x_{2}|\ \ \ \ \mathrm{for} \ \mathrm{any} \ \ x_{1}, x_{2} \in [0, \Lambda/ \mu],$$
where $g(x)=G(x)/ x$.
\end{enumerate}

Notice that the function $G(I)$ in this paper is not necessarily increasing function with respect to $I$,
e.g., it may firstly increase and then decrease
describing some psychological effects: for a very large number of
infectives the infection force $\beta G(I)$ may decrease as $I$ increases,
because in the presence of a very large number of infectives the population
may tend to reduce the number of contacts per unit time \cite{Capasso1978}.
Our general results can be used in
some specific forms for the incidence rate that have been commonly used, for example:

(\textrm{i}) linear type: $G(I)=I$;

(\textrm{ii}) saturated incidence rate: $G(I)= I/(1+ a I)$ (e.g., \cite{Capasso1978});

(\textrm{iii}) non-monotonic incidence rate: $G(I)= I/(1+ a I^{2})$ (e.g., \cite{Ruan2007});

(\textrm{iv}) incidence rates with ``media coverage'' effect as shown below:

{\bf Type 1} (see \cite{Cui1}): $G(I)=I \exp(-mI)$, where $m$ is a positive constant.

{\bf Type 2} (see \cite{Cui2}):
 $\beta G(I)=(\beta-\widetilde{\beta} f(I))I$,
where $\beta>\widetilde{\beta}$ and the twice-differentiable function $f(I)$ satisfies
$f(0)=0,\ f^{'}(I)\geq 0,\ \lim_{I\rightarrow +\infty}f(I)=1$,
where $f'(0)$ denotes the derivative of the function $f(x)$ at $x=0$.

To study the dynamics of model (\ref{2.2}), we need the following two lemmas,
the proof of which are straightforward, so are omitted.

\begin{lem}
\label{lemma-global}
Suppose that the system
$$\frac{d x(t)}{dt}=F(x(t)),\ \ t\geq 0 $$
with $F: \mathbb{R}^{3}\rightarrow \mathbb{R}^{3}$
has a globally asymptotically stable equilibrium $x^{*}\in \mathbb{R}^{3}$.
Then, for any neighborhood $U$ of $x^{*}$ and any compact set $C\subset \mathbb{R}^{3}$,
there exists a $T>0$ such that $x(t, x_{0})\in U$ for all $t\geq T$ and $x_{0}\in C$.
\end{lem}

\begin{lem}
\label{lemma2.1.} For any given initial
value $\big(z_{0}, r(0)\big) \in \mathbb{R}^{3}_{+}\times \mathcal{M}$,
there exists a unique solution $(S(t, z_{0})$, $I(t, z_{0})$, $R(t, z_{0}))$
of system (\ref{2.2}) on $t\geq 0$ and the solution will
always remain in $\mathbb{R}^{3}_{+}$.
Moreover, the solution has the property
that for every $\omega \in \widetilde{\Omega}$,
there exists $t_{0}=t_{0}(\omega)$ such that
\begin{eqnarray*}
 \frac{\Lambda}{\mu+\alpha} < S(t,\omega,z_{0})+I(t,\omega,z_{0})+R(t,\omega,z_{0})
 <\frac{\Lambda}{\mu}, \ \ \mathrm{for} \ \mathrm{all}\ t\geq t_{0}.
\end{eqnarray*}
\end{lem}

Let
\begin{eqnarray}\label{state-space}
   X= \mathcal{K}\times \mathcal{M},
\end{eqnarray}
where
$ \mathcal{K}=\big\{x\in \mathbb{R}^{3}_{+}: \Lambda/(\mu+\alpha)<x_{1}+x_{2}+x_{3}< \Lambda/\mu \big\}$.
From Lemma \ref{lemma2.1.}, we get that
the set $X$
is an attraction domain for the system (\ref{2.2}) in the sense that
all sample paths of system (\ref{2.2}) tend to this set
and once they go into it, they will remain there forever.
Hence, without loss of generality, we consider
only the smaller state-space $X$ throughout this paper.

\subsection{Main results}

\subsubsection{Extinction and persistence of the disease}

In studying epidemic modelings,
the most interesting and important
issues are usually to establish the threshold condition for the extinction and persistence
of the disease,
which will be given in this subsection.
Let us introduce the basic reproduction number in a random environment
as follows
\begin{eqnarray*}\label{R0}
\mathcal{R}_{0}=\frac{\sum_{e\in \mathcal{M}}\pi_{e}(\Lambda \beta_{e}G'(0)/\mu)}{
\mu+\alpha+\delta}.
\end{eqnarray*}
Here, we refer the readers to \cite{Gray2012,Bacar2013,Bacar2014} for the reason why we label the
above formula as the basic reproduction number of our stochastic epidemic model (\ref{2.2}).
Next, we shall show that the position of $\mathcal{R}_{0}$ with respect to $1$ serves as a threshold
between the disease extinction
and persistence for our epidemic model (\ref{2.2}) of SIRS type in a random environment.

Before we prove the results, let us state a proposition
which gives an equivalent condition for the position of $\mathcal{R}_{0}$ with respect to $1$ in terms
of the system parameters and the stationary distribution $\pi$ of the Markov chain $r(t)$, which
drives the switch of environments.

\begin{prop} \label{proposition3.1.}
The following alternative conditions on the value of $\mathcal{R}_{0}$ are valid:

(i) $\mathcal{R}_{0}<1$ if and only if $\sum_{e\in \mathcal{M}}\pi_{e}B(e)<0$;

(ii) $\mathcal{R}_{0}>1$ if and only if $\sum_{e\in \mathcal{M}}\pi_{e}B(e)>0$,\\
where
$$B(e)=\frac{\Lambda \beta_{e}G'(0)}{\mu} - (\mu+\alpha+\delta)$$
for each $e \in \mathcal{M}$.
\end{prop}

The proof of this proposition is straightforward, so is omitted.
We then state our results on the extinction of the disease.

\begin{thm}\label{theorem3.1.}
Suppose that $\mathcal{R}_{0}<1$, then the solution $(S(t), I(t), R(t))$
of system (\ref{2.2}) with any initial value
$\big(z_{0}, r(0)\big)\in X$ has the property that
\begin{eqnarray}
\lim_{t\rightarrow +\infty}S(t)&=& \Lambda/\mu \ \ \ \ \ \ \ a.s., \label{3.1}\\ [+5pt]
\lim_{t\rightarrow +\infty}I(t)&=& 0\ \ \ \ \ \ \ a.s., \label{3.2}\\  [+5pt]
\lim_{t\rightarrow +\infty}R(t)&=& 0\ \ \ \ \ \ \ a.s. \label{3.3}
\end{eqnarray}
\end{thm}

\begin{pf} From the second equation of system (\ref{2.2}),
it is easy to see that
\begin{eqnarray*}
\frac{ d\log I(t)}{dt}&=& \beta_{r(t)}\frac{SG(I)}{I} -\big(\mu+\alpha+\delta \big).
\end{eqnarray*}
By Lemma \ref{lemma2.1.} and the assumption (\textbf{H1}): $G(I)\leq I G'(0)$,
we have
\begin{eqnarray*}
\frac{ d\log I(t)}{dt}&\leq& \frac{\Lambda \beta_{r(t)} G'(0)}{\mu} -\big(\mu+\alpha+\delta \big).
\end{eqnarray*}
Integrating the above inequality, it is obtained from the Birkhoff Ergodic theorem
that
\begin{eqnarray*}
 \limsup_{t\rightarrow +\infty}\frac{\log I(t)}{t}\leq
 \sum_{e\in \mathcal{M}}\pi_{e}B(e)\ \ \ \ a.s.,
\end{eqnarray*}
this, together with (i) of Proposition \ref{proposition3.1.}, implies
\begin{eqnarray*}
 \lim_{t\rightarrow +\infty}I(t)&=& 0\ \ \ \ \ \ \ a.s.
\end{eqnarray*}
This is the required assertion (\ref{3.2}).
The procedure to prove assertions (\ref{3.1}) and (\ref{3.3})
is similar to that given in Theorem 1 in \cite{Dan2015},
so is omitted. This completes the proof of Theorem \ref{theorem3.1.}.\ \ $\Box$
\end{pf}

\begin{rmk}\label{remark3.1.}
From the proof procedure of Theorem \ref{theorem3.1.},
it is obvious to see that when $\mathcal{R}_{0}<1$, any positive solution of system (\ref{2.2})
converges exponentially to
the disease-free state $(\Lambda/\mu, 0,0)$ with probability $1$.
\end{rmk}

We now turn to the persistence of the disease.
\begin{thm}\label{theorem3.2.}
Suppose that $\mathcal{R}_{0}>1$,
then for any initial value $\big(z_{0}, r(0)\big)\in X$,
the following statement is valid with probability $1$:
\begin{eqnarray*}
  \liminf_{t\rightarrow +\infty} \frac{1}{t}\int_{0}^{t}I(s) ds \geq
         \frac{\mu^{2}}{\beta^{M}(\mu \vartheta+\beta^{M}(G'(0))^{2})\Lambda}
         \sum_{e\in \mathcal{M}}\pi_{e}B(e).
\end{eqnarray*}
By (ii) of Proposition \ref{proposition3.1.}, we hence conclude that the disease is persistent in the time mean with probability $1$.
\end{thm}

The proof of this theorem is similar to the arguments of Theorem 2 in \cite{Dan2015},
so we omit it.

\begin{rmk}
\label{Remark3.2.} Indeed, from Theorem \ref{theorem3.2.}, it can be seen that
the persistence in the time mean implies the following weak persistence.
That is, if $\mathcal{R}_{0}>1$,
then for any initial value $\big(z_{0}, r(0)\big)\in X$,
\begin{eqnarray*}
\limsup_{t\rightarrow +\infty} I(t)\geq
         \frac{\mu^{2}}{\beta^{M}(\mu \vartheta+\beta^{M}(G'(0))^{2})\Lambda}
         \sum_{e\in \mathcal{M}}\pi_{e}B(e)>0 \ \ \ \ a.s.
\end{eqnarray*}
\end{rmk}

Let us now establish a useful corollary, which indicates that
$\mathcal{R}_{0}$ is a threshold value determining the disease is extinct or persistent, i.e.,
the position of the deterministic quantity $\mathcal{R}_{0}$
with respect to $1$ determines
the disease extinction or persistence  for system (\ref{2.2}).
We can easily obtain from Theorems \ref{theorem3.1.} and \ref{theorem3.2.} the following corollary.

\begin{cor}
\label{corollary3.1.}
For any initial value $\big(z_{0}, r(0)\big)\in X$,
the solution $(S(t), I(t), R(t))$
of system (\ref{2.2}) has the property that

(\textrm{i}) if $\mathcal{R}_{0}<1$, the number of infected individuals $I(t)$ of system (\ref{2.2})
tends to zero exponentially almost surely, i.e., the disease dies out with
probability one;

(\textrm{ii}) if $\mathcal{R}_{0}>1$, the disease will be almost surely persistent
in the time mean.
\end{cor}

\subsubsection{$\Omega$-limit set and attractor}

For each state $e\in \mathcal{M}$,
we denote by $\pi_{t}^{e}(z_{0})$ the solution of
system (\ref{2.2}) in the state $e$ with the initial value $z_{0}\in \mathcal{K}$.
As in \cite{Yingang2011}, the $\Omega$-limit set of the trajectory starting from
an initial value $z_{0}\in \mathcal{K}$ is defined by
$$\Omega(z_{0}, \omega)=\bigcap_{T>0} \overline{\bigcup_{t>T}(S(t,\omega, z_{0}),I(t,\omega, z_{0}),R(t,\omega, z_{0}))}.$$
We here use the notation ``$\Omega$-limit set'' in lieu of the usual one ``$\omega$-limit set'' in
the deterministic dynamical system for
avoiding notational conflict with the element notation $\omega$ in the probability sample space.
In this subsection,
we shall show that under some appropriate conditions,
$\Omega(z_{0}, \omega)$ is deterministic, i.e., it is constant almost surely;
moreover, it is independent of the initial value $z_{0}$.

Let us recall that the basic reproduction number $\mathcal{R}_{0}^{e}$ of the deterministic subsystem of (\ref{2.2})
corresponding to environmental state $e$
can be computed by using the next generation matrix approach \cite{Brauer2008} as
\begin{eqnarray}\label{D-Rsub}
\mathcal{R}_{0}^{e}=\frac{\Lambda\beta_{e} G'(0)}{\mu(\mu+\alpha+\delta)}.
\end{eqnarray}
In the remainder of this paper, we sometimes need to impose the following assumption:
\begin{itemize}
   \item [(\textbf{H3})] If $\mathcal{R}_{0}^{e}> 1$, then there exists a unique and
    globally asymptotically stable positive
   equilibrium $E^{*}_{e} = (S^{*}_{e}, I^{*}_{e},R^{*}_{e})$
   for the corresponding deterministic subsystem of (\ref{2.2}).
\end{itemize}

Note that the condition $\mathcal{R}_{0}>1$
implies that there exists at least one state $e\in \mathcal{M}$ such that
$B(e)>0$, i.e., $\mathcal{R}_{0}^{e}>1$ corresponding to the subsystem of (\ref{2.2})
in the state $e$. We may assume, without any loss of generality,
that $\mathcal{R}_{0}^{1}>1$ in the remainder of this paper
if the assumption (\textbf{H3}) is needed.
Indeed, Remark \ref{remark3} indicates that
the assumption is reasonable.
Hence, the subsystem in the first state ($e=1$) has a globally
stable positive equilibrium $E^{*}_{1} = (S^{*}_{1}, I^{*}_{1},R^{*}_{1})$.

Now we recall some concepts on the Lie algebra of vector fields \cite{Jurdjevic1997,Rudnicki-Pichor}.
Let $a(x)$ and $b(x)$ be two vector fields on $\mathbb{R}^{d}$. The
Lie bracket $[a, b]$ is also a vector field given by
$$[a, b]_{j}(x)=\sum_{k=1}^{d}\Bigg(a_{k}\frac{\partial b_{j}}{\partial x_{k}}(x)
 -b_{k}\frac{\partial a_{j}}{\partial x_{k}}(x)\Bigg),\ \ j=1,2,\ldots, d.$$
In the remainder of this paper, we need the following definition:
A point $z=(S,I,R)\in \mathbb{R}^{3}_{+}$ is said to satisfy the condition ($\mathbf{H}$),
if vectors $Y_{1}(z), \ldots, Y_{E}(z)$, $[Y_{i}$, $Y_{j}](z)_{i,j \in \mathcal{M}}$,
$[Y_{i}$, $[Y_{j}$, $Y_{k}]](z)_{i,j,k \in \mathcal{M}}$, $\ldots$,
span the space $\mathbb{R}^{3}$, where for each $e\in \mathcal{M}$,
$$Y_{e}(S,I,R)=\left(
          \begin{array}{c}
            \Lambda-\mu S+\lambda R-\beta_{e} S G(I) \\
            \beta_{e} S G(I)-\big(\mu + \alpha +\delta \big)I \\
            \delta I-\big(\mu +\lambda\big)R \\
          \end{array}
        \right).
$$
With this definition, we have
\begin{thm}\label{theorem4.1.}
Suppose that $\mathcal{R}_{0}>1$ and the hypotheses (\textbf{H3}) hold.
Let
$$\Gamma=\Bigg\{(S,I,R)=\pi_{t_{k}}^{p_{k}}\circ \cdots \circ \pi_{t_{1}}^{p_{1}}(E^{*}_{1}): t_{1}, \ldots, t_{k}\geq 0 \ \mathrm{and} \ p_{1}, \ldots, p_{k}\in \mathcal{M}, k\in \mathbb{N} \Bigg \}.$$
Then, the following statements are valid:

(a) With $\overline{\Gamma}$ denoting the closure of $\Gamma$,
$\overline{\Gamma}$ is a subset of the $\Omega$-limit set $\Omega(z_{0}, \omega)$
with probability $1$.

(b) If there exists a point $\overline{z}_{*}:=(S_{*}, I_{*}, R_{*})\in \Gamma$ satisfying
the condition (\textbf{H}), then $\Gamma$ absorbs all positive solutions
in the sense that for any initial value $z_{0}\in \mathcal{K}$,
the value
$$\widetilde{T}(\omega)=\inf \Big \{t>0: (S(s,\omega,z_{0}), I(s,\omega,z_{0}), R(s,\omega,z_{0}))\in \Gamma, \ \forall s>t \Big \}$$
is finite outside a $\mathbb{P}$-null set. Consequently,
$\overline{\Gamma}$ is the $\Omega$-limit set $\Omega(z_{0}, \omega)$
for any $z_{0}\in \mathcal{K}$ with probability $1$.
\end{thm}

\begin{rmk}\label{remark3}
If there exists the other state $e\in \mathcal{M}$
such that $\mathcal{R}_{0}^{e}>1$, we can define
$$\Gamma_{e}=\Bigg\{(S,I,R)=\pi_{t_{k}}^{p_{k}}\circ \cdots \circ \pi_{t_{1}}^{p_{1}}(E^{*}_{e}): t_{1}, \ldots, t_{k}\geq 0 \ \mathrm{and} \ p_{1}, \ldots, p_{k}\in \mathcal{M}, k\in \mathbb{N} \Bigg \}.$$
Then, we have the same conclusion that
the closure $\overline{\Gamma}_{e}$ of $\Gamma_{e}$
is the $\Omega$-limit set $\Omega(z_{0}, \omega)$
with probability $1$, which implies that $\overline{\Gamma}_{e}=\overline{\Gamma}$.
\end{rmk}

\subsubsection{Global convergence of the distribution in total variation norm}
\label{Stability in the distribution}

First, we establish the existence of an invariant probability measure.
\begin{thm}\label{theorem5.1.}
If $\mathcal{R}_{0}>1$, the Markov process $((S(t), I(t), R(t)), r(t))$
has an invariant probability measure $\nu^{*}$
on the state space $X$.
\end{thm}

\begin{rmk}
To prove this theorem, we utilize the alternative principle in \cite{Meyn-III,Stettner1986}.
In the proceeding of the proof,
it is critical to find a compact subset $O$ of $\overline{X}$
such that the probability that the process $((S(t), I(t), R(t)), r(t))$
falls into $O$ in the time mean is positive,
where $\overline{X}$ is a slightly larger state space
compared with the state space $X$.
In Section \ref{Proofs of main results},
we shall give two alternative methods for constructing
such compact subset: one is similar to that
given in Theorem 3.1 of \cite{Yingang2014},
while the other one is based on the results
of Theorem \ref{theorem3.2.}.
\end{rmk}

We now characterize the invariant probability
measure by the following theorem.

\begin{thm}\label{theorem5.2.}
Suppose that $\mathcal{R}_{0}>1$ and the hypotheses (\textbf{H3}) hold.
Assume further that the hypotheses in (b) of Theorem \ref{theorem4.1.} is satisfied.
Then for any initial value
$x_{0}=\big(z_{0}, r(0)\big)\in X$,
the instant distribution of the process $((S(t), I(t), R(t)), r(t))$ satisfies
\begin{eqnarray}\label{theorem19-1}
\lim_{t\rightarrow +\infty} \|\mathbb{P}(t, x_{0}, (\cdot, \cdot))-\nu^{*}(\cdot, \cdot) \| &=& 0,
\end{eqnarray}
and
\begin{eqnarray}\label{theorem19-2}
\mathbb{P}_{x_{0}}\Bigg\{\lim_{t\rightarrow +\infty} \frac{1}{t}
\int_{0}^{t}f(z(s,z_{0}), r(s))ds=\sum_{e\in\mathcal{M}}\int_{\mathcal{K}}
f(u,e)\nu^{*}(du,e) \Bigg\} &=& 1,
\end{eqnarray}
where $f(\cdot, \cdot)$ is any $\nu^{*}$-integrable function.
Moreover,
the stationary distribution $\nu^{*}$ has the density $f^{*}$
with respect to the product measure $m$ on $X$
and $\mathrm{supp}(f^{*})= \Gamma\times \mathcal{M}$.

\end{thm}

\begin{rmk}
In \cite{Yingang2014}, it is easy to see that if the condition (4.4) is valid,
then the condition (\textbf{H}) holds.
In (b) of Theorem  \ref{theorem4.1.} and Theorem \ref{theorem5.2.},
we therefore have weakened the conditions ensuring
the global attractivity of the $\Omega$-limit set of the system
and the convergence in total variation of the
instantaneous measure to the stationary measure. Moreover,
the method in \cite{Yingang2011,Yingang2014} is not
applicable to the case where the dimension of the system considered
is higher than the number of environmental
regimes, while the modified method in proving Theorems \ref{theorem4.1.} and \ref{theorem5.2.}
has no limitations
for both the number of environmental
regimes and the dimension of the system considered.
This, together with the following example \ref{exam}, illustrates that our modified method
in this paper has a wider application than the techniques in \cite{Yingang2011,Yingang2014}.
\end{rmk}

\begin{figure}[htbp]
\begin{center}
  \includegraphics[height=3.6 in, width=4.8 in]{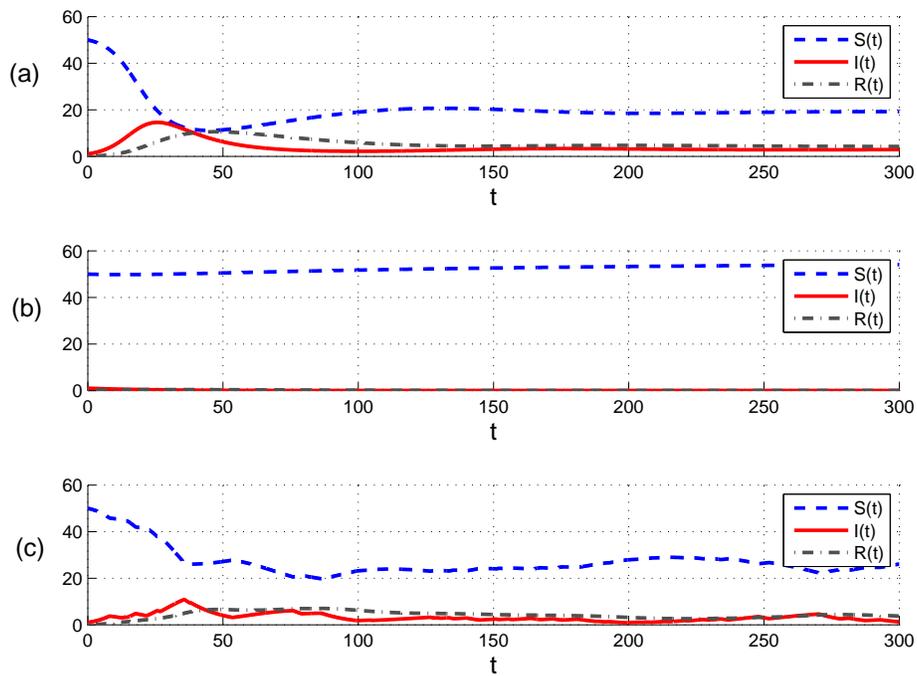}
  \caption{The paths of $S(t)$, $I(t)$ and $R(t)$ for the systems
  (\ref{example01}) and (\ref{example02}) with the same initial values $z_{0}=(50, 1, 0)$.
  (a) the paths of the deterministic system (\ref{example01}) in the state 1
  with $\beta_{1}=0.0056\  \mathrm{days}^{-1}$ corresponding to $\mathcal{R}^{1}_{0}=2.9057>1$;
  (b) the paths of the deterministic system (\ref{example01}) in the state 2
  with $\beta_{2}=0.0013\  \mathrm{days}^{-1}$ corresponding to $\mathcal{R}^{2}_{0}=0.6745<1$;
  (c) the paths of the stochastic system (\ref{example02}) switching between the state 1 and the state 2 with $\mathcal{R}_{0}=1.8726 >1$ and $r(0)=1$.
  (Color figure online)}\label{Fig1}
  \end{center}
\end{figure}

\begin{exam}
\label{exam}

Let $G(I)=I/(1+a I^{2})$, where
$a$ is a positive parameter measuring the psychological or inhibitory effect.
The stochastic model
(\ref{2.2}) under regime switching reduces to
\begin{eqnarray}\label{example02}
\left\{
  \begin{array}{l}
d S(t)/dt=\Lambda-\mu S(t)+\lambda R(t)-\beta_{r(t)} S(t) \frac{I(t)}{1+a I^{2}(t)},   \\ [+8pt]
d I(t)/dt=\beta_{r(t)} S(t) \frac{I(t)}{1+a I^{2}(t)}-\big(\mu + \alpha +\delta \big)I(t), \\ [+8pt]
d R(t)/dt=\delta I(t)-\big(\mu +\lambda\big)R(t),
\end{array} \right.
\end{eqnarray}
which switches from one to the other according to the movement
of the right-continuous Markov chain
$\{r(t), t\geq0 \}$
taking values in the state space $\mathcal{M}=\{1, 2\}$
with the generator
$$Q=\frac{1-0.5}{365}\cdot \left(
      \begin{array}{cc}
        -169 & 169 \\
        196 & -196 \\
      \end{array}
    \right).
$$
This Markov chain has a unique stationary distribution
\begin{eqnarray}\label{distri}
\pi &=& (\pi_{1}, \pi_{2})=\Big( \frac{196}{365}, \frac{169}{365}\Big).
\end{eqnarray}

For each environmental state $e\in \mathcal{M}$,
Cai and Kang et al. \cite{Cai2015}
discussed the corresponding deterministic epidemic model:
\begin{eqnarray}\label{example01}
\left\{
  \begin{array}{l}
d S(t)/dt=\Lambda-\mu S(t)+\lambda R(t)-\beta_{e} S(t) \frac{I(t)}{1+a I^{2}(t)},   \\ [+8pt]
d I(t)/dt=\beta_{e} S(t) \frac{I(t)}{1+a I^{2}(t)}-\big(\mu + \alpha +\delta \big)I(t), \\ [+8pt]
d R(t)/dt=\delta I(t)-\big(\mu +\lambda\big)R(t).
\end{array} \right.
\end{eqnarray}
The basic reproduction number of the deterministic system (\ref{example01}) is
$$\mathcal{R}_{0}^{e}=\frac{\Lambda\beta_{e} G'(0)}{\mu(\mu+\alpha+\delta)}=\frac{\Lambda \beta_{e}}{\mu(\mu + \alpha +\delta)},$$
which determines the extinction and persistence of the disease.
According to Theorem 5.1 in \cite{Cai2015}, then
\begin{itemize}
    \item[$\bullet$] The unique disease-free
    equilibrium $E_{0} = (\Lambda/\mu, 0, 0)$ is globally asymptotically stable
    whenever $\mathcal{R}^{e}_{0} \leq 1$, and it is unstable when $\mathcal{R}^{e}_{0} > 1$.
    \item[$\bullet$] When $\mathcal{R}^{e}_{0}> 1$, there exists a unique
    equilibrium $E^{*}_{e} = (S^{*}_{e}, I^{*}_{e},R^{*}_{e})$ with $S_{e}^{*}>0$, $I_{e}^{*}>0$, $R_{e}^{*}>0$, which is
    globally asymptotically stable.
\end{itemize}
Without any loss of generality, let us assume
that $\mathcal{R}_{0}^{1}>1$ and $\mathcal{R}_{0}^{2}<1$ in this example.

The other parameter values used here are mainly taken from the work in \cite{Anderson1979}
investigating the dynamics of Pasteurella muris in colonies of laboratory
mice:
$\Lambda=0.33 \ \mathrm{days}^{-1}$, $\mu=0.006\  \mathrm{days}^{-1}$,
$\alpha=0.06\  \mathrm{days}^{-1}$, $\delta=0.04\  \mathrm{days}^{-1}$,
$\lambda=0.021\  \mathrm{days}^{-1}$, $a=0.001$ \cite{Cai2015},
$\beta_{1}=0.0056\  \mathrm{days}^{-1}$ corresponding to $\mathcal{R}^{1}_{0}=2.9057>1$
(the disease is persistent, see (a) in Figure \ref{Fig1}), and
$\beta_{2}=0.0013\  \mathrm{days}^{-1}$ corresponding to $\mathcal{R}^{2}_{0}=0.6745<1$
(the disease is extinct, see (b) in Figure \ref{Fig1}).
Combining with (\ref{distri}), this implies that
the basic reproduction number for the stochastic system (\ref{example02}) is
\begin{eqnarray*}
\mathcal{R}_{0}=\frac{\pi_{1}\Lambda \beta_{1}}{\mu (
\mu+\alpha+\delta)}+\frac{\pi_{2}\Lambda \beta_{2}}{\mu (
\mu+\alpha+\delta)}=1.8726 >1.
\end{eqnarray*}
Hence, by Corollary \ref{corollary3.1.}, the disease will be almost surely persistent
in the time mean (see (c) in Figure \ref{Fig1}).
For the system (\ref{example02}), we define the following set as in Theorem \ref{theorem4.1.}
$$\Gamma=\Bigg\{(S,I,R)=\pi_{t_{k}}^{p_{k}}\circ \cdots \circ \pi_{t_{1}}^{p_{1}}(E^{*}_{1}): t_{1}, \ldots, t_{k}\geq 0 \ \mathrm{and} \ p_{1}, \ldots, p_{k}\in \{1, 2\}, k\in \mathbb{N} \Bigg \},$$
where $E^{*}_{1}=(19.0161, 2.8783, 4.2830)$.

For any point $z=(S,I,R)\in \Gamma$, let
$$Y_{e}(z)=\left(
          \begin{array}{c}
            \Lambda-\mu S+\lambda R-\beta_{e} S G(I) \\
            \beta_{e} S G(I)-\big(\mu + \alpha +\delta \big)I \\
            \delta I-\big(\mu +\lambda\big)R \\
          \end{array}
        \right), \ \ \ e=1,\ 2.
$$
Because there are only two environmental states for the system (\ref{example02})
with three equations,
the matrix $\big(Y_{1}(z), Y_{2}(z) \big)$ has $3$ rows and $2$ columns,
which indicates that the determinant
$$\mathrm{det}\Big(Y_{1}(z), Y_{2}(z) \Big)$$
becomes meaningless.
Hence, the method used in \cite{Yingang2014}
can not be used to judge whether or not the claims that
the $\Omega$-limit set $\Gamma$ of the system (\ref{example02}) is globally attractive
and its instantaneous measure converges in total
variation to  some stationary measure are valid,
while our method in this paper can be feasible.

Computing the Lie bracket of $Y_{1}(z)$ and $Y_{2}(z)$, we obtain the vector field given by
$$[Y_{1}, Y_{2}](z)=(\beta_{1}-\beta_{2})\left(
                      \begin{array}{c}
                        (\Lambda+\lambda R)G(I)-(\mu+\alpha+\delta)SIG'(I) \\
                        -(\Lambda+\lambda R)G(I)-(\alpha+\delta)SG(I)+(\mu+\alpha+\delta)SIG'(I) \\
                        \delta SG(I) \\
                      \end{array}
                    \right).
$$
Hence,
\begin{eqnarray*}
 && \mathrm{det}\Big(Y_{1}(z), Y_{2}(z), [Y_{1}, Y_{2}](z) \Big)\nonumber \\[+6pt]
 &=& -\Big[(\beta_{1}-\beta_{2})SG(I)\Big]^{2}\cdot \Big[
\mu \delta \Big( \frac{\Lambda}{\mu}-(S+I+R)\Big)-\alpha (\mu+\lambda)R\Big]\nonumber \\[+6pt]
 &=& -\Big[\frac{(\beta_{1}-\beta_{2})SI}{1+a I^{2}}\Big]^{2}\cdot \Big[
\mu \delta \Big( \frac{\Lambda}{\mu}-(S+I+R)\Big)-\alpha (\mu+\lambda)R\Big],
\end{eqnarray*}
from which, we can easily find by the numerical method that
there exist many points $z\in \Gamma$
such that $\mathrm{det}\big(Y_{1}(z), Y_{2}(z), [Y_{1}, Y_{2}](z) \big)\neq 0$.
For instance, $\mathrm{det}\big(Y_{1}(z), Y_{2}(z), [Y_{1}, Y_{2}](z) \big)\neq 0$
when $z=\pi_{100}^{2}(E^{*}_{1})=(37.3966, 0.0033, 0.4464)$, which implies that
at the point $z=\pi_{100}^{2}(E^{*}_{1})$, the vectors
$Y_{1}(z)$, $Y_{2}(z)$, $[Y_{1}, Y_{2}](z)$
span the space $\mathbb{R}^{3}$, i.e.,
the point $z=\pi_{100}^{2}(E^{*}_{1})$ satisfies the condition (\textbf{H}).
Therefore, by (b) of Theorem  \ref{theorem4.1.} and Theorem \ref{theorem5.2.},
we can conclude that
the $\Omega$-limit set $\Gamma$ of the system (\ref{example02}) is globally attractive
and its instantaneous measure converges to some stationary measure
in total variation (see, Figures \ref{Fig2} and \ref{Fig3}).

\end{exam}

\begin{figure}[htbp]
\begin{center}
  \includegraphics[height=3.2 in, width=3.8
  in]{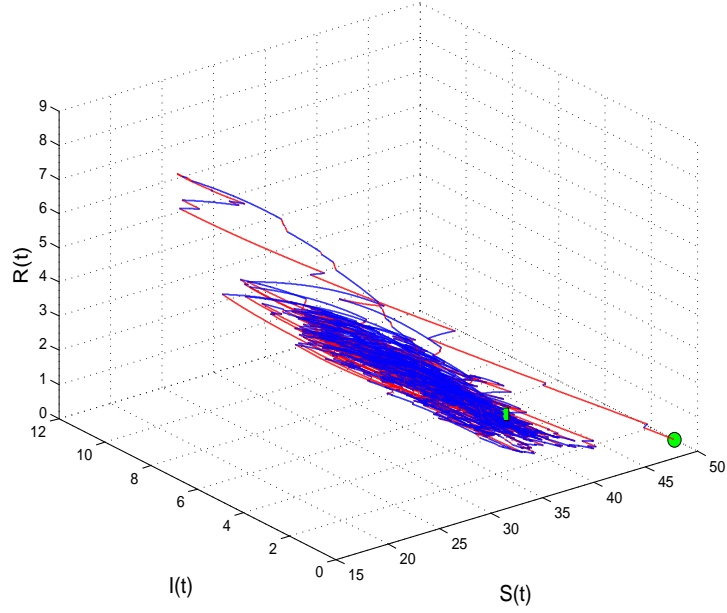}
  \caption{A sample orbit of the system (\ref{example02})
  with the initial values $z_{0}=(50, 1, 0)$ and $r(0)=1$. The red lines and blue
  lines represent the paths of the system (\ref{example02})
  in states $1$ and $2$, respectively. The $\circ$ denotes the
  the starting point of the orbit, the $\square$ denotes the end
  of the orbit, and the number of environmental switching is $2000$.
  (Color figure online)}\label{Fig2}
  \end{center}
\end{figure}

\begin{figure}[htbp]
\begin{center}
  \includegraphics[height=3.6 in, width=4.8
  in]{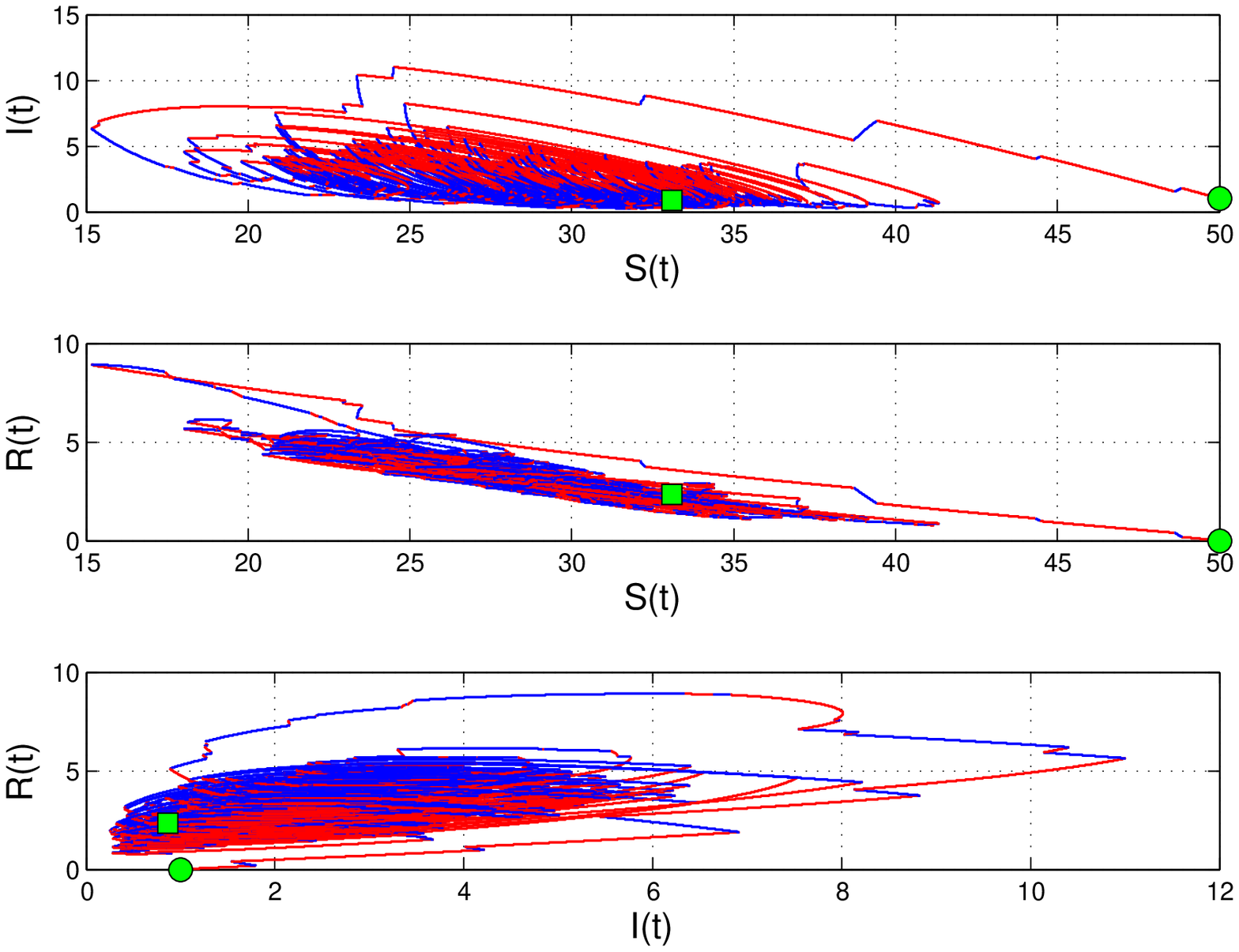}
  \caption{Two-dimension projections of the sample orbit in Figure \ref{Fig2}.
  The orbit of Figure \ref{Fig2} is
  projected onto the S-I, S-R and I-R coordinate planes from the top to the bottom, respectively.
  (Color figure online)}\label{Fig3}
  \end{center}
\end{figure}

\section{Proofs of main results}
\label{Proofs of main results}

In this section, we first provide some auxiliary definitions and results concerning
stability of Markovian processes \cite{Meyn-I,Meyn-II,Meyn-III,Stettner1986}
that we will use later to prove our main results.

Let $\mathbf{\Phi}=\{ \mathbf{\Phi}_{t}: t\geq 0\}$
be a time-homogeneous Markov process with state space $(\mathbf{X}, \mathcal{B}(\mathbf{X}))$,
where the state space $\mathbf{X}$ is a locally compact and separable metric space, and
$\mathcal{B}(\mathbf{X})$ is the Borel $\sigma$-algebra on $\mathbf{X}$.
The process $\mathbf{\Phi}$ evolves on the probability space
$(\widetilde{\Omega}, \mathcal{F}, \mathbb{P}_{x})$,
where $x\in \mathbf{X}$ is the initial condition of the process, and the measure
$\mathbb{P}_{x}$ satisfies $\mathbb{P}_{x}(\mathbf{\Phi}_{t}\in B)=\mathbb{P}(t, x, B)$
for all $x\in \mathbf{X}$, $t\geq0$ and $B\in \mathcal{B}(\mathbf{X})$.
Assume further that $\mathbf{\Phi}$ is a Borel right process,
so that in particular $\mathbf{\Phi}$ is a strongly Markovian process
with right-continuous sample paths \cite{SHARPE1988}. Note that if $\mathbf{\Phi}$
is a Feller process then it is a Borel right process.
Let $\{\mathbb{P}^{t}\}_{t\geq 0}$ be the transition semigroup of the process $\mathbf{\Phi}$,
and the operator $\mathbb{P}^{t}$ acts on $\sigma$-finite
measure $\nu$ on $\mathbf{X}$ via
$$\nu \mathbb{P}^{t}(B)=\int_{\mathbf{X}}\mathbb{P}(t, x, B)\nu(dx),\ \ B\in \mathcal{B}(\mathbf{X}).$$
It is easy to see that if $\nu$ is the initial distribution of the process $\mathbf{\Phi}$ at time $0$,
then $\nu \mathbb{P}^{t}$ is the distribution of $\mathbf{\Phi}$ at time $t$.
A $\sigma$-finite measure $\nu$ on $\mathcal{B}(\mathbf{X})$, with the property
$$\nu(B)=\nu \mathbb{P}^{t}(B)=\int_{\mathbf{X}}\mathbb{P}(t, x, B)\nu(dx)$$
for any $B\in \mathcal{B}(\mathbf{X})$ and $t\geq0$,
will be called invariant (or stationary) measure.

For a measurable set $B\in \mathcal{B}(\mathbf{X})$, let
$$\tau_{B}=\inf \{t\geq0: \mathbf{\Phi}_{t}\in B\}, \ \ \ \ \eta_{B}=\int_{0}^{\infty}\mathbf{1}_{\{\mathbf{\Phi}_{t}\in B\}}dt.$$
We then introduce the standard definition of Harris recurrence as follows.
The process $\mathbf{\Phi}$ is called Harris recurrent if either

(a) some nontrivial $\sigma$-finite
measure $\phi_{1}$ exists such that $\mathbb{P}_{x}\{\tau_{B}<\infty\}\equiv1$ whenever $\phi_{1}(B)>0$; or

(b) some nontrivial $\sigma$-finite
measure $\phi_{2}$ exists such that $\mathbb{P}_{x}\{\eta_{B}=\infty\}\equiv1$ whenever $\phi_{2}(B)>0$, \\
is valid for any $x\in \mathbf{X}$ and $B\in \mathcal{B}(\mathbf{X})$.
In addition, the process $\mathbf{\Phi}$ is said to be positive Harris recurrent
if it is Harris recurrent with a finite invariant measure.

Suppose that $\mathfrak{a}$ is a probability measure on $[0, +\infty)$,
and define the transition function $\mathbf{K}_{\mathfrak{a}}$ of a general
sampled Markov chain corresponding to
the measure $\mathfrak{a}$ as
$$\mathbf{K}_{\mathfrak{a}}(x,B)=\int_{0}^{\infty}\mathbb{P}(t,x,B)\mathfrak{a}(dt)$$
for all $x\in \mathbf{X}$ and $B\in \mathcal{B}(\mathbf{X})$.
We refer the readers to \cite{Meyn-II} for the detail definition
of the sampled Markov chain of the process $\mathbf{\Phi}$.
A kernel $\mathcal{T}: \mathbf{X}\times \mathcal{B}(\mathbf{X})\rightarrow [0, +\infty)$
is called a continuous component of the Markov transition function $\mathbf{K}_{\mathfrak{a}}$
if

(\textrm{i}) For each $B\in \mathcal{B}(\mathbf{X})$ the function $\mathcal{T}(\cdot, B)$ is lower semi-continuous;

(\textrm{ii}) For all $x\in \mathbf{X}$ and $B\in \mathcal{B}(\mathbf{X})$,
the measure $\mathcal{T}(x, \cdot)$ satisfies $\mathbf{K}_{\mathfrak{a}}(x, B) \geq \mathcal{T}(x, B)$.\\
The continuous component $\mathcal{T}$ is called everywhere non-trivial if $\mathcal{T}(x, \mathbf{X}) > 0$
for each $x\in \mathbf{X}$.
The process $\mathbf{\Phi}$ will be called a $\mathcal{T}$-process if there exists a probability measure $\mathfrak{a}$, such that the corresponding
transition function $\mathbf{K}_{\mathfrak{a}}$ has an everywhere non-trivial continuous component $\mathcal{T}$.

The process $\mathbf{\Phi}$ is called bounded in probability on average
if for each initial condition $x\in \mathbf{X}$ and any $\varepsilon>0$,
there exists a compact subset $C \subset \mathbf{X}$ such that
$$\liminf_{t\rightarrow +\infty} \frac{1}{t}\int_{0}^{t}\mathbb{P}_{x}
\{\mathbf{\Phi}_{s}\in C\}ds\geq 1-\varepsilon.$$

For the $\sigma$-finite measure $\phi$,
the process $\mathbf{\Phi}$ is called $\phi$-irreducible if
for any $x\in \mathbf{X}$ and $B\in \mathcal{B}(\mathbf{X})$,
$$\mathbb{E}_{x}[\eta_{B}]>0$$
whenever $\phi(B)>0$, where $\mathbb{E}_{x}(\cdot)$
denotes the expectation of the random variable with respect to the probability
measure $\mathbb{P}_{x}$. In this case, the measure $\phi$ is called
an irreducibility measure. The process $\mathbf{\Phi}$
is usually be called irreducible when
the specific irreducibility measure is irrelevant.

For checking the asymptotic stability of the distribution, we need some results concerning
the stochastic stability and ergodicity of Markovian processes as follows.

\begin{thm}\label{invariant-exist} (see, \cite{Meyn-III,Stettner1986})
Assume that $\mathbf{\Phi}$ is a Feller process. Then either

(a) there exists an invariant probability measure on $\mathbf{X}$, or

(b) for any compact set $C \subset \mathbf{X}$,
$$\lim_{t\rightarrow \infty} \sup_{\nu} \frac{1}{t}\int_{0}^{t}
\Bigg( \int_{\mathbf{X}}\mathbb{P}(s, x, C)\nu(dx)\Bigg)ds=0,$$
where the supremum is taken over all initial distributions $\nu$ on the state space $\mathbf{X}$.
\end{thm}

\begin{thm}\label{theorem-stability-1} (see, Theorem 3.2 in \cite{Meyn-II})
Suppose that the process $\mathbf{\Phi}$ is an irreducible $\mathcal{T}$-process.
Then $\mathbf{\Phi}$ is positive Harris recurrent if and only if
$\mathbf{\Phi}$ is bounded in probability on average.

\end{thm}

\begin{thm}\label{theorem-stability-2} (see, Theorem 8.1 in \cite{Meyn-II})
Suppose that the process $\mathbf{\Phi}$ is irreducible, and is
bounded in probability on average. Assume further
that there exists a lattice distribution $\mathfrak{a}$, such that
the corresponding kernel $\mathbf{K}_{\mathfrak{a}}$ possesses an
everywhere non-trivial continuous component. Then the following
assertions are valid.

(i) For every $x\in \mathbf{X}$,
$$\lim_{t\rightarrow +\infty} \|\mathbb{P}(t, x, \cdot)-\Pi(\cdot) \|=0,$$
where $\|\cdot\|$ denotes the total variation norm,
and $\Pi$ denotes the invariant probability measure on the state space $\mathbf{X}$.

(ii) Moreover, for any $\Pi$-integrable function $f$ and any $x\in \mathbf{X}$,
$$\mathbb{P}_{x}\Bigg\{\lim_{t\rightarrow +\infty} \frac{1}{t}
\int_{0}^{t}f(\mathbf{\Phi}_{s})ds=\int_{\mathbf{X}}f(u)\Pi(du) \Bigg\}=1.$$

\end{thm}

From the preceding definitions, it is easy to get that
if the process $\mathbf{\Phi}$ is Harris recurrent,
then it is also irreducible. Hence, the following proposition is immediate
from Theorems \ref{theorem-stability-1} and \ref{theorem-stability-2}.

\begin{prop}
\label{corollary5.1.}
Suppose that the process $\mathbf{\Phi}$ is positive Harris recurrent.
If a lattice distribution $\mathfrak{a}$ exists such that
the corresponding kernel $\mathbf{K}_{\mathfrak{a}}$ admits an
everywhere non-trivial continuous component, i.e.,
$\mathbf{\Phi}$ is a $\mathcal{T}$-process,
then we have

(i) For every $x\in \mathbf{X}$,
$$\lim_{t\rightarrow +\infty}\|\mathbb{P}(t, x, \cdot)-\Pi(\cdot) \|=0,$$
where $\Pi$ denotes the invariant probability measure on $\mathbf{X}$.

(ii) Moreover, for any $\Pi$-integrable function $f$ and any $x\in \mathbf{X}$,
$$\mathbb{P}_{x}\Bigg\{\lim_{t\rightarrow +\infty} \frac{1}{t}
\int_{0}^{t}f(\mathbf{\Phi}_{s})ds=\int_{\mathbf{X}}f(u)\Pi(du) \Bigg\}=1.$$

\end{prop}

Let $\mathcal{B}(X)$ be the $\sigma$-algebra
of Borel subsets of the space $X$ defined by (\ref{state-space}),
and $\widetilde{m}$ be the Lebesgue measure on $\mathcal{K}$;
$\ell$ be the measure on $\mathcal{M}$ given by $\ell(e)=\pi_{e}$, $e\in \mathcal{M}$.
Denote by $m=\widetilde{m}\times \ell$ be the product measure on $(X,\mathcal{B}(X))$
given by $m(B\times \{e \})=\widetilde{m}(B)\times \ell(e)$ for each
$B\in \mathcal{B}(\mathcal{K})$ and $e \in \mathcal{M}$.

\subsection{Proof of Theorem \ref{theorem4.1.}}

To prove Theorem \ref{theorem4.1.}, we need the following three lemmas.

\begin{lem}\label{lemma4.1.}
Suppose that $\mathcal{R}_{0}>1$.
For any initial condition $(z_{0}, r(0))\in X$,
there exists a $\iota> 0$ such that
$\limsup_{t\rightarrow +\infty} I(t, z_{0})\geq \iota$ \ \ a.s.
\end{lem}

\begin{pf}
By Remark \ref{Remark3.2.}, the result of Lemma \ref{lemma4.1.}
is immediate. However, we then give another proving method.
Let
$$\varepsilon=\frac{\mu \sum_{e\in \mathcal{M}}\pi_{e}B(e)}{4\Lambda \beta^{M}}.$$
By the assumptions (\textbf{H1}) and (\textbf{H2}),
we get that $\lim_{x\rightarrow 0^{+}}G(x)/x =G'(0)$, which
yields that there exists a positive constant $\iota$ satisfying
$$\iota<\frac{\mu^{2} \sum_{e\in \mathcal{M}}\pi_{e}B(e)}{8\Lambda (\beta^{M}G'(0))^{2}},$$
such that $G'(0)-G(x)/x< \varepsilon$ if $0<x<\iota$.
Next, we shall show that $\limsup_{t\rightarrow +\infty} I(t, z_{0})> \iota$ \ a.s.
In the contrary, assume that a measurable set $B$ with $\mathbb{P}(B)>0$
exists such that $\limsup_{t\rightarrow +\infty} I(t,\omega, z_{0})< \iota$
for any $\omega\in B$. Then, for each $\omega\in B$, there exists
a $T=T(\omega, \iota)>0$ such that $I(t,\omega, z_{0})< \iota$
for all $t>T$. Note that $G(I(t,\omega, z_{0}))\leq I(t,\omega, z_{0})G'(0)$
and $S(t,\omega, z_{0})<\Lambda/\mu$ for all $t>T$ and $\omega\in B$,
which implies that
$$\beta_{r(t,\omega)}S(t,\omega, z_{0})G(I(t,\omega, z_{0}))\leq \frac{\Lambda \beta^{M}G'(0)\iota}{\mu}$$
for all $t>T$ and $\omega\in B$.
Hence, it is obtained from the first equation of system (\ref{2.2})
that for each $\omega\in B$,
\begin{eqnarray*}
\frac{dS(t,\omega, z_{0})}{dt} &=& \Lambda-\mu S(t,\omega, z_{0})
+\lambda R(t,\omega, z_{0})- \beta_{r(t,\omega)}S(t,\omega, z_{0})G(I(t,\omega, z_{0}))\\[+6pt]
&\geq& \Lambda-\mu S(t,\omega, z_{0})-\frac{\Lambda \beta^{M}G'(0)\iota}{\mu}
\end{eqnarray*}
holds for all $t>T$. By the comparison theorem, one can find a $T_{1}=T_{1}(\omega,\iota)>T$ satisfying
\begin{eqnarray}\label{4.1}
 S(t,\omega, z_{0})&\geq& \frac{\Lambda}{\mu}- \frac{2\Lambda \beta^{M}G'(0)\iota}{\mu^{2}}
\end{eqnarray}
for all $t>T_{1}$ and $\omega\in B$. Noting that
$G'(0)-(G(I(t,\omega, z_{0}))/I(t,\omega, z_{0}))<\varepsilon$ for all $t>T_{1}$ and $\omega\in B$,
it follows from (\ref{4.1}) that
\begin{eqnarray*}
&&\beta_{r(t,\omega)}\Bigg( \frac{\Lambda G'(0)}{\mu}
-\frac{S(t,\omega, z_{0})G(I(t,\omega, z_{0}))}{I(t,\omega, z_{0})}\Bigg) \\[+6pt]
&=& \beta_{r(t,\omega)}\Bigg( \frac{\Lambda G'(0)}{\mu}
- S(t,\omega, z_{0})G'(0)\Bigg)+\beta_{r(t,\omega)}\Bigg(S(t,\omega, z_{0})G'(0)\\[+6pt]
&&-\frac{S(t,\omega, z_{0})G(I(t,\omega, z_{0}))}{I(t,\omega, z_{0})}\Bigg) \\[+6pt]
&\leq & \beta^{M}G'(0)\Big(\frac{\Lambda}{\mu}-S(t,\omega, z_{0}) \Big)
+\frac{\Lambda \beta^{M}}{\mu}\Bigg(G'(0)- \frac{G(I(t,\omega, z_{0}))}{I(t,\omega, z_{0})}\Bigg) \\[+6pt]
&\leq& \frac{\sum_{e\in \mathcal{M}}\pi_{e}B(e)}{4}+ \frac{\sum_{e\in \mathcal{M}}\pi_{e}B(e)}{4}\\[+6pt]
&=&\frac{\sum_{e\in \mathcal{M}}\pi_{e}B(e)}{2}
\end{eqnarray*}
for all $t>T_{1}$ and $\omega\in B$.
From the second equation of system (\ref{2.2}),
we then have
\begin{eqnarray*}
d \ln I(t,\omega, z_{0})
&=& \Bigg[ \beta_{r(t,\omega)}\frac{S(t,\omega, z_{0})G(I(t,\omega, z_{0}))}{I(t,\omega, z_{0})}-
(\mu+\alpha+\delta) \Bigg]dt \\[+6pt]
&=& \Bigg(\frac{\Lambda\beta_{r(t,\omega)} G'(0)}{\mu}-
(\mu+\alpha+\delta) \Bigg)dt\\[+6pt]
&&- \beta_{r(t,\omega)}\Bigg(\frac{\Lambda G'(0)}{\mu}-\frac{S(t,\omega, z_{0})G(I(t,\omega, z_{0}))}{I(t,\omega, z_{0})} \Bigg)dt  \\[+6pt]
&\geq& \Big(\frac{\Lambda\beta_{r(t,\omega)} G'(0)}{\mu}-
(\mu+\alpha+\delta) \Big)dt-\frac{\sum_{e\in \mathcal{M}}\pi_{e}B(e)}{2}dt
\end{eqnarray*}
for all $t>T_{1}$ and $\omega\in B$.
Integrating both sides of the above inequality from $T_{1}$ to $t$ ($t>T_{1}$) yields
\begin{eqnarray}\label{4.2}
&&\frac{\ln I(t,\omega, z(T_{1}, \omega, z_{0}))-\ln I(T_{1}, \omega, z_{0})}{t} \nonumber \\[+6pt]
&\geq& \frac{1}{t} \int_{T_{1}}^{t}\Big(\frac{\Lambda\beta_{r(s,\omega)} G'(0)}{\mu}-
(\mu+\alpha+\delta) \Big)ds-\frac{\sum_{e\in \mathcal{M}}\pi_{e}B(e)}{2}
\end{eqnarray}
for all $t>T_{1}$ and $\omega\in B$,
where $z(T_{1}, \omega, z_{0})=(S(T_{1}, \omega, z_{0}),I(T_{1}, \omega, z_{0}),R(T_{1}, \omega, z_{0}))$.
Since $I(t,\omega, z(T_{1}, \omega, z_{0}))$ is bounded,
$$\limsup_{t\rightarrow +\infty}\frac{\ln I(t,\omega, z(T_{1}, \omega, z_{0}))-\ln I(T_{1}, \omega, z_{0})}{t}\leq 0.$$
However, by the Birkhoff Ergodic theorem, it is obtained from the right hand of
the inequality (\ref{4.2}) that
$$\lim_{t\rightarrow +\infty}\Bigg[ \frac{1}{t} \int_{T_{1}}^{t}\Big(\frac{\Lambda\beta_{r(s,\omega)} G'(0)}{\mu}-
(\mu+\alpha+\delta) \Big)ds-\frac{\sum_{e\in \mathcal{M}}\pi_{e}B(e)}{2} \Bigg]=\frac{\sum_{e\in \mathcal{M}}\pi_{e}B(e)}{2}>0$$
for almost all $\omega\in B$.
This is a contradiction. Thus, $\limsup_{t\rightarrow +\infty} I(t, z_{0})> \iota$ \ a.s.
This completes the proof of Lemma \ref{lemma4.1.}.\ \ $\Box$

\end{pf}

\begin{lem}\label{lemma4.2.}
If $\mathcal{R}_{0}>1$,
then there exists a $\varrho>0$ such that
for each $i\in \mathcal{M}$, the event
$$\mathbf{E}_{i} = \Big\{\mathrm{for}\  \mathrm{infinite} \
\mathrm{many}\  k\in \mathbb{N}_{0}, (S(\tau_{k}),I(\tau_{k}),R(\tau_{k}))\in H_{\varrho}\
\mathrm{occur} \ \mathrm{with}\ r(\tau_{k})=i \Big\}$$
occurs with probability one, where
$$H_{\varrho}=\Big\{(S,I,R):\frac{\Lambda}{\mu+\alpha}\leq
S+I+R\leq \frac{\Lambda}{\mu},\ I\geq \varrho \Big\}.$$
\end{lem}

\begin{pf}
Fix a $T>0$, by Lemma \ref{lemma4.1.}, one can define
almost surely finite stopping times with respect
to filtration $\mathcal{F}_{t}$:
\begin{eqnarray*}
\eta_{1}&=& \inf \Big\{t>0: I(t)\geq \iota \Big\},  \\[+6pt]
\eta_{2}&=& \inf \Big\{t>\eta_{1}+T: I(t)\geq \iota \Big\}, \\[+6pt]
&\cdots&  \\[+6pt]
\eta_{n}&=& \inf \Big\{t>\eta_{n-1}+T: I(t)\geq \iota \Big\},\ \ldots.
\end{eqnarray*}
Let $\tau(\eta_{k})=\inf \{t>\eta_{k}:r(t)\neq r(\eta_{k})\}$,
$\widehat{\sigma}(\eta_{k})=\tau(\eta_{k})-\eta_{k}$
and $A_{k}=\{\widehat{\sigma}(\eta_{k})<T\}$, $k\in \mathbb{N}$.
Then $A_{k+1}$ is in the $\sigma$-algebra generated by
$\{r(\eta_{k+1}+s):s\geq 0 \}$ while $A_{k}\in \mathcal{F}_{\eta_{k+1}}$.
From the strong Markov property of the process $((S(t),I(t),R(t)), r(t))$,
it is obtained that for any $A_{k}^{c}$, $k\in \mathbb{N}$,
\begin{eqnarray*}
\mathbb{P}\Big(A_{k}^{c}\Big) &=& \sum_{e\in \mathcal{M}}\mathbb{P}\Big(A_{k}^{c}\ \Big| \ r(\eta_{k})=e\Big)\mathbb{P}\Big(r(\eta_{k})=e\Big)  \\[+2pt]
&=& \sum_{e\in \mathcal{M}}\mathbb{P}\Big(\widehat{\sigma}(\eta_{k})\geq T\ \Big|\  r(\eta_{k})=e\Big)\mathbb{P}\Big(r(\eta_{k})=e\Big)  \\[+2pt]
&=&\sum_{e\in \mathcal{M}}\mathbb{P}\Big(\widehat{\sigma}(0)\geq T\ \Big| \ r(0)=e\Big)\mathbb{P}\Big(r(\eta_{k})=e\Big)  \\[+2pt]
&\leq& \overline{p}_{1}
\end{eqnarray*}
and
\begin{eqnarray*}
 \mathbb{P}\Big(A_{k}^{c}\cap A_{k+1}^{c}\Big)&=& \sum_{e\in \mathcal{M}}\mathbb{P}\Big(A_{k}^{c}\cap A_{k+1}^{c}\ \Big| \ r(\eta_{k+1})=e \Big)\mathbb{P}\Big(r(\eta_{k+1})=e\Big)  \\[+2pt]
&=& \sum_{e\in \mathcal{M}}\mathbb{P}\Big(A_{k}^{c}\ \Big|\  r(\eta_{k+1})=e \Big)\mathbb{P}\Big(A_{k+1}^{c}\ \Big|\  r(\eta_{k+1})=e\Big)\mathbb{P}\Big(r(\eta_{k+1})=e\Big)  \\[+2pt]
&=& \sum_{e\in \mathcal{M}}\mathbb{P}\Big(A_{k}^{c}\ \Big| \ r(\eta_{k+1})=e\Big)\mathbb{P}\Big(\widehat{\sigma}(\eta_{k+1})\geq T\ \Big|\  r(\eta_{k+1})=e\Big)\mathbb{P}\Big(r(\eta_{k+1})=e\Big)  \\[+2pt]
&\leq& \overline{p}_{1} \sum_{e\in \mathcal{M}}\mathbb{P}\Big(A_{k}^{c}\ \Big|\  r(\eta_{k+1})=e\Big)\mathbb{P}\Big(r(\eta_{k+1})=e\Big) \\[+2pt]
&=&\overline{p}_{1}^{2},
\end{eqnarray*}
where $\overline{p}_{1}:=\max_{e\in \mathcal{M}}\{\mathbb{P}(\widehat{\sigma}(0)\geq T| r(0)=e)\}<1$.
This implies
$$\mathbb{P}\Bigg(\bigcap_{i=1}^{\infty}\bigcup_{k=i}^{\infty}A_{k} \Bigg)=
1-\mathbb{P}\Bigg(\bigcup_{i=1}^{\infty}\bigcap_{k=i}^{\infty}A_{k}^{c} \Bigg)=1,$$
that is, the events $A_{k}$, $k\in \mathbb{N}$ occur infinitely often a.s.
Since
\begin{eqnarray}\label{fangsuo}
\frac{dI(t)}{dt} &=& \beta_{r(t)}S(t)G(I(t))-(\mu+\alpha+\delta)I(t)\geq -(\mu+\alpha+\delta)I(t)
\end{eqnarray}
with $I(\eta_{k})\geq \iota$,
it follows from the comparison theorem
that $I(\eta_{k}+t)\geq \iota e^{-(\mu+\alpha+\delta)t}$ for all $t\geq 0$.
Thus, with probability one the events $\{(S(\eta_{k}+\widehat{\sigma}(\eta_{k})), I(\eta_{k}+\widehat{\sigma}(\eta_{k})), R(\eta_{k}+\widehat{\sigma}(\eta_{k})))\in H_{\widetilde{\varrho}}\}$,
$k\in \mathbb{N}$ occur infinitely often with $\widetilde{\varrho}:=\iota e^{-(\mu+\alpha+\delta)T}$,
which yields that
the events
$\{(S(\tau_{k}), I(\tau_{k}), R(\tau_{k}))\in H_{\widetilde{\varrho}}\}$,
$k\in \mathbb{N}_{0}$ occur infinitely often a.s.

Let us fix $i\in \mathcal{M}$. Introduce the following sequence
of stopping times with respect to
filtration $\mathcal{F}_{0}^{n}$:
\begin{eqnarray*}
\xi_{1}&=& \inf \Big\{k\in \mathbb{N}_{0}:\  (S(\tau_{k}), I(\tau_{k}), R(\tau_{k}))\in H_{\widetilde{\varrho}} \Big\}, \\[+6pt]
\xi_{2}&=& \inf \Big\{k\in \mathbb{N}_{0}:\  \tau_{k}>\tau_{\xi_{1}}+T, (S(\tau_{k}), I(\tau_{k}), R(\tau_{k}))\in H_{\widetilde{\varrho}} \Big\}, \\[+6pt]
&\cdots&  \\[+6pt]
\xi_{n}&=& \inf \Big\{k\in \mathbb{N}_{0}:\ \tau_{k}>\tau_{\xi_{n-1}}+T, (S(\tau_{k}), I(\tau_{k}), R(\tau_{k}))\in H_{\widetilde{\varrho}} \Big\},\ \ldots.
\end{eqnarray*}
Define the events
$$\tilde{A}_{k}=\Big \{\mathrm{for} \ \mathrm{some}\  s\in (0,T],
r((\tau_{\xi_{k}}+s)-)\neq r(\tau_{\xi_{k}}+s)\
\mathrm{and}\  r(\tau_{\xi_{k}}+s)=i\Big \},$$
$k\in \mathbb{N}$.
Next, to obtain the assertion that $\mathbb{P}\{\mathbf{E}_{i}\}=1$,
by (\ref{fangsuo}) we need only to prove that with probability one
the events $\tilde{A}_{k}$, $k\in \mathbb{N}$ occur infinitely often
provided we set $\varrho:=\widetilde{\varrho} e^{-(\mu+\alpha+\delta)T}$.
From the strong Markov property of the process $((S(t),I(t),R(t)), r(t))$,
it then follows that for any $\tilde{A}_{k}^{c}$, $k\in \mathbb{N}$,
\begin{eqnarray*}
\mathbb{P}\Big(\tilde{A}_{k}^{c}\Big)
&=& \sum_{n=0}^{\infty}\mathbb{P}\Big(\tilde{A}_{k}^{c}\ \Big| \ \xi_{k}=n\Big)\mathbb{P}\Big(\xi_{k}=n\Big)  \\[+2pt]
&=& \sum_{n=0}^{\infty}\mathbb{P}\Big(\xi_{k}=n\Big)\cdot
\sum_{e\in \mathcal{M}}\mathbb{P}\Big(\tilde{A}_{k}^{c}\ \Big| \ \xi_{k}=n, r(\tau_{\xi_{k}})=e\Big)\mathbb{P}\Big(r(\tau_{\xi_{k}})=e\ \Big| \ \xi_{k}=n\Big)  \\[+2pt]
&=& \sum_{n=0}^{\infty}\mathbb{P}\Big(\xi_{k}=n\Big)\cdot
\sum_{e\in \mathcal{M}}\mathbb{P}\Big(\tilde{A}_{k}^{c}\ \Big| \ \xi_{k}=n, r(\tau_{n})=e\Big)\mathbb{P}\Big(r(\tau_{\xi_{k}})=e\ \Big| \ \xi_{k}=n\Big)  \\[+2pt]
&=& \sum_{n=0}^{\infty}\mathbb{P}\Big(\xi_{k}=n\Big)\cdot
\sum_{e\in \mathcal{M}}\mathbb{P}\Big(\mathbf{E}\ \Big| \ r(0)=e\Big)\mathbb{P}\Big(r(\tau_{\xi_{k}})=e\ \Big| \ \xi_{k}=n\Big)  \\[+2pt]
&\leq& \widetilde{p}_{1}
\sum_{n=0}^{\infty}\mathbb{P}\Big(\xi_{k}=n\Big)\cdot
\sum_{e\in \mathcal{M}}\mathbb{P}\Big(r(\tau_{\xi_{k}})=e\ \Big| \ \xi_{k}=n\Big)  \\[+2pt]
&=&\widetilde{p}_{1}
\end{eqnarray*}
and
\begin{eqnarray*}
 &&\mathbb{P}\Big(\tilde{A}_{k}^{c}\cap \tilde{A}_{k+1}^{c}\Big)\\[+2pt]
 &=& \sum_{n=0}^{\infty}\mathbb{P}\Big(\tilde{A}_{k}^{c}\cap \tilde{A}_{k+1}^{c}\ \Big| \ \xi_{k+1}=n \Big)\mathbb{P}\Big(\xi_{k+1}=n\Big)  \\[+2pt]
 &=& \sum_{n=0}^{\infty}\mathbb{P}\Big(\xi_{k+1}=n\Big)\cdot
\sum_{e\in \mathcal{M}}\mathbb{P}\Big(\tilde{A}_{k}^{c}\cap \tilde{A}_{k+1}^{c}\ \Big| \ \xi_{k+1}=n, r(\tau_{\xi_{k+1}})=e\Big)\mathbb{P}\Big(r(\tau_{\xi_{k+1}})=e\ \Big| \ \xi_{k+1}=n\Big)  \\[+2pt]
&=& \sum_{n=0}^{\infty}\mathbb{P}\Big(\xi_{k+1}=n\Big)\cdot
\sum_{e\in \mathcal{M}}\mathbb{P}\Big(\tilde{A}_{k}^{c}\ \Big| \ \xi_{k+1}=n, r(\tau_{\xi_{k+1}})=e\Big)
\mathbb{P}\Big(\tilde{A}_{k+1}^{c}\ \Big| \ \xi_{k+1}=n, r(\tau_{\xi_{k+1}})=e\Big)\\[+2pt]
&&\cdot \mathbb{P}\Big(r(\tau_{\xi_{k+1}})=e\ \Big| \ \xi_{k+1}=n\Big)  \\[+2pt]
&=& \sum_{n=0}^{\infty}\mathbb{P}\Big(\xi_{k+1}=n\Big)\cdot
\sum_{e\in \mathcal{M}}\mathbb{P}\Big(\tilde{A}_{k}^{c}\ \Big| \ \xi_{k+1}=n, r(\tau_{\xi_{k+1}})=e\Big)
\mathbb{P}\Big(\mathbf{E}\ \Big| \ r(0)=e\Big)\\[+2pt]
&&\cdot \mathbb{P}\Big(r(\tau_{\xi_{k+1}})=e\ \Big| \ \xi_{k+1}=n\Big)  \\[+2pt]
&\leq& \widetilde{p}_{1} \sum_{n=0}^{\infty}\mathbb{P}\Big(\xi_{k+1}=n\Big)\cdot
\sum_{e\in \mathcal{M}}\mathbb{P}\Big(\tilde{A}_{k}^{c}\ \Big| \ \xi_{k+1}=n, r(\tau_{\xi_{k+1}})=e\Big)
\mathbb{P}\Big(r(\tau_{\xi_{k+1}})=e\ \Big| \ \xi_{k+1}=n\Big)\\[+2pt]
&\leq& \widetilde{p}_{1} \sum_{n=0}^{\infty}\mathbb{P}\Big(\tilde{A}_{k}^{c}\ \Big| \ \xi_{k+1}=n \Big)\mathbb{P}\Big(\xi_{k+1}=n\Big)\\[+2pt]
&\leq&\widetilde{p}_{1}^{2},
\end{eqnarray*}
where $\widetilde{p}_{1}:=\max_{e\in \mathcal{M}}\{\mathbb{P}(\mathbf{E}\ | \ r(0)=e)\}<1$
with the event $$\mathbf{E}=\Big \{\mathrm{for} \ \mathrm{all}\  s\in (0,T], \mathrm{if} \
r(s-)\neq r(s)\
\mathrm{then}\  r(s)\neq i\Big \}.$$
Hence, by induction we have
$$\mathbb{P}\Bigg(\bigcap_{i=1}^{\infty}\bigcup_{k=i}^{\infty}\tilde{A}_{k} \Bigg)=
1-\mathbb{P}\Bigg(\bigcup_{i=1}^{\infty}\bigcap_{k=i}^{\infty}\tilde{A}_{k}^{c} \Bigg)=1,$$
which means that the events $\tilde{A}_{k}$, $k\in \mathbb{N}$ occur infinitely often a.s.
The proof of Lemma \ref{lemma4.2.} is completed.\ \ $\Box$

\end{pf}

\begin{lem}\label{lemma4.3.}
Fix any two distinct states $e_{1}, e_{2}\in \mathcal{M}$ satisfying
the probability that the state $e_{2}$ can be accessible from
the state $e_{1}$ through only one step is positive, i.e., $q_{e_{1},e_{2}}>0$.
Let $\{\xi_{n}\}_{n=1}^{\infty}$ be a sequence of strictly increasing finite stopping times with respect to filtration $\mathcal{F}_{0}^{n}$ satisfying $r(\tau_{\xi_{n}})=e_{1}$ for all $n\in \mathbb{N}$.
Suppose that $B$ is a bounded Borel subset of
$[0, +\infty)$
with positive Lebesgue measure.
Then, the events $A_{k}=\{r(\tau_{\xi_{k}+1})=e_{2}, \ \sigma_{\xi_{k}+1}\in B\}, k\in \mathbb{N}$ occur infinitely often a.s., that is,
$$\mathbb{P}\Bigg(\bigcap_{i=1}^{\infty}\bigcup_{k=i}^{\infty}A_{k} \Bigg)=1.$$
\end{lem}

\begin{pf}
By the strong Markov property of the process
$((S(t),I(t),R(t)), r(t))$,
it is obtained that for each $k\in \mathbb{N}$,
\begin{eqnarray*}
\mathbb{P}\Big(A_{k}\Big) &=& \sum_{n=0}^{\infty}\mathbb{P}\Big(A_{k}\ \Big| \ \xi_{k}=n\Big)\mathbb{P}\Big(\xi_{k}=n\Big)  \\[+2pt]
&=& \sum_{n=0}^{\infty}\mathbb{P}\Big(\xi_{k}=n\Big)\cdot
\sum_{e\in \mathcal{M}}\mathbb{P}\Big(A_{k}\ \Big| \ \xi_{k}=n, r(\tau_{\xi_{k}})=e\Big)\mathbb{P}\Big(r(\tau_{\xi_{k}})=e\ \Big| \ \xi_{k}=n\Big)  \\[+2pt]
&=& \sum_{n=0}^{\infty}\mathbb{P}\Big(\xi_{k}=n\Big)
\mathbb{P}\Big(A_{k}\ \Big| \ \xi_{k}=n, r(\tau_{n})=e_{1}\Big)\mathbb{P}\Big(r(\tau_{\xi_{k}})=e_{1}\ \Big| \ \xi_{k}=n\Big)  \\[+2pt]
&=& \sum_{n=0}^{\infty}\mathbb{P}\Big(\xi_{k}=n\Big)
\mathbb{P}\Big(r(\tau_{n+1})=e_{2}, \ \sigma_{n+1}\in B\ \Big| \ r(\tau_{n})=e_{1}\Big)\\[+2pt]
&=& \sum_{n=0}^{\infty}\mathbb{P}\Big(\xi_{k}=n\Big)
\mathbb{P}\Big(r(\tau_{1})=e_{2}, \ \sigma_{1}\in B\ \Big| \ r(0)=e_{1}\Big)\\[+2pt]
&=&
\mathbb{P}\Big(r(\tau_{1})=e_{2}, \ \sigma_{1}\in B\ \Big| \ r(0)=e_{1}\Big):=\widetilde{p}>0.
\end{eqnarray*}
This implies that $\mathbb{P}(A_{k}^{c})=1-\mathbb{P}(A_{k})
=1-\widetilde{p}<1$. Moreover,
\begin{eqnarray*}
&&\mathbb{P}\Big(A_{k}^{c}\cap A_{k+1}^{c}\Big)\\[+2pt]
&=& \sum_{n=0}^{\infty}\mathbb{P}\Big(A_{k}^{c}\cap A_{k+1}^{c}\ \Big| \ \xi_{k+1}=n \Big)\mathbb{P}\Big(\xi_{k+1}=n\Big)  \\[+2pt]
&=& \sum_{n=0}^{\infty}\mathbb{P}\Big(\xi_{k+1}=n\Big)\cdot
\sum_{e\in \mathcal{M}}\mathbb{P}\Big(A_{k}^{c}\cap A_{k+1}^{c}\ \Big| \ \xi_{k+1}=n, r(\tau_{\xi_{k+1}})=e\Big)\mathbb{P}\Big(r(\tau_{\xi_{k+1}})=e\ \Big| \ \xi_{k+1}=n\Big)  \\[+2pt]
&=& \sum_{n=0}^{\infty}\mathbb{P}\Big(\xi_{k+1}=n\Big)
\mathbb{P}\Big(A_{k}^{c}\cap A_{k+1}^{c}\ \Big| \ \xi_{k+1}=n, r(\tau_{\xi_{k+1}})=e_{1}\Big)\mathbb{P}\Big(r(\tau_{\xi_{k+1}})=e_{1}\ \Big| \ \xi_{k+1}=n\Big)  \\[+2pt]
&=& \sum_{n=0}^{\infty}\mathbb{P}\Big(\xi_{k+1}=n\Big)
\mathbb{P}\Big(A_{k}^{c}\ \Big| \ \xi_{k+1}=n, r(\tau_{\xi_{k+1}})=e_{1}\Big)\mathbb{P}\Big(A_{k+1}^{c}\ \Big| \ \xi_{k+1}=n, r(\tau_{\xi_{k+1}})=e_{1}\Big)  \\[+2pt]
&=& (1-\widetilde{p})\sum_{n=0}^{\infty}\mathbb{P}\Big(\xi_{k+1}=n\Big)
\mathbb{P}\Big(A_{k}^{c}\ \Big| \ \xi_{k+1}=n, r(\tau_{\xi_{k+1}})=e_{1}\Big)\mathbb{P}\Big(r(\tau_{\xi_{k+1}})=e_{1}\ \Big| \ \xi_{k+1}=n\Big)  \\[+2pt]
&=& (1-\widetilde{p})\sum_{n=0}^{\infty}\mathbb{P}\Big(\xi_{k+1}=n\Big)
\mathbb{P}\Big(A_{k}^{c}\ \Big| \ \xi_{k+1}=n\Big)  \\[+2pt]
&=& (1-\widetilde{p})^{2}
\end{eqnarray*}
for all $k\in \mathbb{N}$. Hence, by induction we have
$$\mathbb{P}\Bigg(\bigcup_{i=1}^{\infty}\bigcap_{k=i}^{\infty}A_{k}^{c} \Bigg)=0.$$
Thus,
$$\mathbb{P}\Bigg(\bigcap_{i=1}^{\infty}\bigcup_{k=i}^{\infty}A_{k} \Bigg)=1-
\mathbb{P}\Bigg(\bigcup_{i=1}^{\infty}\bigcap_{k=i}^{\infty}A_{k}^{c} \Bigg)=1.$$
This completes the proof of Lemma \ref{lemma4.3.}.\ \ $\Box$

\end{pf}

\begin{pot1}

(a) Let us divide the following proof of the assertion (a)
into three steps.

{\bf Step 1.} We shall show that $E^{*}_{1}\in \Omega(z_{0}, \omega)$\ a.s.
From Lemma \ref{lemma4.2.}, without loss of generality,
we assume that the event $\mathbf{E}_{1}$ in Lemma \ref{lemma4.2.} occurs.
Specifically, there exists a sequence of strictly increasing finite stopping
times $\{\xi_{n}\}_{n=1}^{\infty}$ with respect to filtration
$\mathcal{F}_{0}^{n}$ satisfying $r(\tau_{\xi_{n}})=1$ and
$(S(\tau_{\xi_{n}}), I(\tau_{\xi_{n}}), R(\tau_{\xi_{n}}))\in H_{\varrho}$
for all $n\in \mathbb{N}$ and almost all $\omega \in \widetilde{\Omega}$.
In addition, for any $\widetilde{\delta}>0$, it is obtained from Lemma \ref{lemma-global}
that one can find a $T>0$ such that $(S(t,z),I(t,z), R(t,z))\in B(E^{*}_{1},\widetilde{\delta})$
is valid for all $t\geq T$ and $z\in \overline{H}_{\varrho}$,
where $\overline{H}_{\varrho}=H_{\varrho} \bigcup \{(S,I,R): S=0, \Lambda/(\mu+\alpha)\leq I+R\leq \Lambda/\mu\} \bigcup \{(S,I,R): R=0, \Lambda/(\mu+\alpha)\leq S+I\leq \Lambda/\mu\}$.
Thus, by Lemma \ref{lemma4.3.},
it follows that with probability one, the events $\{(S(\tau_{\xi_{n}+1}), I(\tau_{\xi_{n}+1}), R(\tau_{\xi_{n}+1}))\in B(E^{*}_{1},\widetilde{\delta})\}$, $k \in \mathbb{N}$ will occur infinitely often with $\sigma_{\xi_{n}+1} \in [T,2T]$.
This means that $E^{*}_{1}\in \Omega(z_{0}, \omega)$\ a.s.

{\bf Step 2.} Let
$$\widetilde{\Gamma}=\Bigg\{(S,I,R)=\pi_{t_{k}}^{p_{k}}\circ \cdots \circ \pi_{t_{1}}^{p_{1}}(E^{*}_{1}): t_{1}, \ldots, t_{k}> 0 \ \mathrm{and} \ p_{1}, \ldots, p_{k}\in \mathcal{M}, k\in \mathbb{N} \Bigg \},$$
with $q_{1,p_{1}}, q_{p_{i},p_{i+1}}>0$ and $p_{i}\neq p_{i+1}$, $i=1,\ldots, k-1$.
We now prove that $\widetilde{\Gamma}\subset \Omega(z_{0}, \omega)$\ a.s.
To begin with, we shall show that for any $t_{1}>0$, $\pi_{t_{1}}^{p_{1}}(E^{*}_{1})\in \Omega(z_{0}, \omega)$\ a.s. Denote $\widetilde{z}^{*}=\pi_{t_{1}}^{p_{1}}(E^{*}_{1})$.
For any $\widetilde{\delta}_{1}>0$, by the continuous dependence of the solutions on
the time and initial conditions, one can find two numbers $\widetilde{\varepsilon}_{1}, \overline{\delta}_{1}>0$
such that if $z\in B(E^{*}_{1}, \widetilde{\varepsilon}_{1})$
then $\pi_{t}^{p_{1}}(z)\in B(\widetilde{z}^{*}, \widetilde{\delta}_{1})$ for all $t\in B(t_{1},\overline{\delta}_{1})$.
Moreover, by Step 1 and Lemma \ref{lemma4.3.}, it follows that
for such number $\widetilde{\varepsilon}_{1}$, there exists a sequence of
strictly increasing finite stopping times $\{\eta_{n}^{0}\}_{n=1}^{\infty}$
with respect to filtration
$\mathcal{F}_{0}^{n}$ satisfying
$$r(\tau_{\eta_{n}^{0}})=p_{1}\  \mathrm{and} \ \Big(S(\tau_{\eta_{n}^{0}}),I(\tau_{\eta_{n}^{0}}),
R(\tau_{\eta_{n}^{0}})\Big)\in B\Big(E^{*}_{1},\widetilde{\varepsilon}_{1}\Big)$$
for all $n\in \mathbb{N}$ and
almost all $\omega\in \widetilde{\Omega}$.
Hence, by Lemma \ref{lemma4.3.} again, we can obtain
a sequence of
strictly increasing finite stopping times $\{\eta_{n}^{1}\}_{n=1}^{\infty}$
with respect to filtration
$\mathcal{F}_{0}^{n}$ satisfying
$$r(\tau_{\eta_{n}^{1}})=p_{2}\  \mathrm{and} \ \Big(S(\tau_{\eta_{n}^{1}}),I(\tau_{\eta_{n}^{1}}),
R(\tau_{\eta_{n}^{1}})\Big)\in B\Big(\widetilde{z}^{*},\widetilde{\delta}_{1}\Big)$$
for all $n\in \mathbb{N}$ and
almost all $\omega\in \widetilde{\Omega}$.
This implies that $\pi_{t_{1}}^{p_{1}}(E^{*}_{1})\in \Omega(z_{0}, \omega)$\ a.s.
Therefore, by induction one can conclude that
$\widetilde{\Gamma}\subset \Omega(z_{0}, \omega)$\ a.s.

{\bf Step 3.} Now, we shall prove that the assertion (a) is valid.
Let us first select any state $\widetilde{e}\in \mathcal{M}$.
If $q_{1,\widetilde{e}}>0$, then by Steps 1 and 2 we have
$\pi_{t}^{\widetilde{e}}(E^{*}_{1})\in \Omega(z_{0}, \omega)$\ a.s.
for all $t\geq 0$. If $q_{1,\widetilde{e}}=0$, we can also
claim that $\pi_{t}^{\widetilde{e}}(E^{*}_{1})\in \Omega(z_{0}, \omega)$\ a.s.
for all $t> 0$. In fact, fix any $t> 0$ and
let $\widetilde{\delta}_{2}$ be a any positive constant.
By the continuous dependence of the solutions on the time and
the initial conditions, one can
find a constant $\widetilde{\varepsilon}_{2}>0$ such that
$\pi_{t}^{\widetilde{e}}(z)\in B(\pi_{t}^{\widetilde{e}}(E^{*}_{1}),\widetilde{\delta}_{2})$
if $z\in B(E^{*}_{1}, \widetilde{\varepsilon}_{2})$.
On the other hand, it is obtained from the irreducibility
of the generator $Q$ that a positive integer $n_{0}$ ($1\leq n_{0} \leq E-1$ )
exists such that
$$q_{l_{0},l_{1}}q_{l_{1},l_{2}}\cdots q_{l_{n_{0} -1},l_{n_{0} }}>0$$
with $\{l_{k}:0\leq k \leq n_{0} \}\subset \mathcal{M}$
and $l_{0}=1$, $l_{n_{0}}=\widetilde{e}$.
Thus, from the continuous dependence of the solutions on the initial conditions,
it follows that there exist $s_{1}, s_{2}, \ldots, s_{n_{0} -1}>0$ such that
$$\pi_{t}^{\widetilde{e}}\circ \pi_{s_{n_{0} -1}}^{l_{n_{0} -1}}\circ \cdots \pi_{s_{1}}^{l_{1}}
(E^{*}_{1})\in B\Big(\pi_{t}^{\widetilde{e}}(E^{*}_{1}), \widetilde{\delta}_{2}\Big).$$
Because $B(\pi_{t}^{\widetilde{e}}(E^{*}_{1}), \widetilde{\delta}_{2})$
is a open set in $\mathbb{R}^{3}_{+}$, a sufficiently small $\widetilde{\delta}_{3}>0$
exists such that
$$B(\overline{z}^{*}, \widetilde{\delta}_{3})
\subset B\Big(\pi_{t}^{\widetilde{e}}(E^{*}_{1}), \widetilde{\delta}_{2}\Big),$$
where $\overline{z}^{*}:= \pi_{t}^{\widetilde{e}}\circ \pi_{s_{n_{0} -1}}^{l_{n_{0} -1}}\circ \cdots \pi_{s_{1}}^{l_{1}}
(E^{*}_{1})$. Hence, it is obtained from Step 2 that
there exists a sequence of
strictly increasing finite stopping times $\{\eta_{n}^{2}\}_{n=1}^{\infty}$
with respect to filtration
$\mathcal{F}_{0}^{n}$,
such that
$$\Big(S(\tau_{\eta_{n}^{2}}),I(\tau_{\eta_{n}^{2}}),
R(\tau_{\eta_{n}^{2}})\Big)\in B(\overline{z}^{*}, \widetilde{\delta}_{3})
\subset B\Big(\pi_{t}^{\widetilde{e}}(E^{*}_{1}), \widetilde{\delta}_{2}\Big)$$
for all $n\in \mathbb{N}$ and
almost all $\omega\in \widetilde{\Omega}$.
This implies that $\pi_{t}^{\widetilde{e}}(E^{*}_{1})\in \Omega(z_{0}, \omega)$\ a.s. for any $t> 0$.
By the similar arguments, we can conclude that
$\Gamma\subset \Omega(z_{0}, \omega)$\ a.s.
Also because $\Omega(z_{0}, \omega)$ is a close set,
we have $\overline{\Gamma}\subset \Omega(z_{0}, \omega)$\ a.s.
This completes the proof of the assertion (a).\\

(b) We shall prove that the second assertion is valid.
By Theorem 1 of Chapter 3 in \cite{Jurdjevic1997},
the claim that the point $\overline{z}_{*}$ satisfies
the condition (\textbf{H}) means that there exist elements
$\hat{p}_{1},\ldots, \hat{p}_{i}$ in $\mathcal{M}$
and $\hat{t}\in \mathbb{R}^{i}$ with positive coordinates
$\hat{t}_{1}, \ldots, \hat{t}_{i}$ such that the function
$\varphi(t_{1},\ldots, t_{i})=\pi_{t_{i}}^{\hat{p}_{i}}\circ \cdots \circ\pi_{t_{1}}^{\hat{p}_{1}}(\overline{z}_{*})$
has a tangent map, the rank of which at $\hat{t}$ is equal
to the dimension of $\mathbb{R}^{3}$.
That is, the rank of matrix
$[\partial \varphi/\partial t_{1}, \ldots, \partial \varphi/\partial t_{i}]
|_{(\hat{t}_{1}, \ldots, \hat{t}_{i})}$
is equal to $3$.
Without loss of generality, we assume that
\begin{eqnarray}\label{4.3}
\mathrm{det} \Big(\frac{\partial \varphi}{\partial t_{1}}, \frac{\partial \varphi}{\partial t_{2}},\frac{\partial \varphi}{\partial t_{3}}\Big)\Big |_{(\hat{t}_{1}, \ldots, \hat{t}_{i})} &\neq& 0.
\end{eqnarray}
By the existence and continuous dependence on the initial conditions of the solutions,
the function
$$\psi(s,t, u)=\pi_{\hat{t}_{i}}^{\hat{p}_{i}}\circ \cdots \circ \pi_{\hat{t}_{4}}^{\hat{p}_{4}}\circ
\pi_{u+\hat{t}_{3}}^{\hat{p}_{3}} \circ \pi_{t+\hat{t}_{2}}^{\hat{p}_{2}}\circ \pi_{s+\hat{t}_{1}}^{\hat{p}_{1}}(\overline{z}_{*})$$
is defined and continuously differentiable in some small domain
$(-a,a)\times (-b,b)\times (-c,c)\subset \mathbb{R}^{3}$.
From (\ref{4.3}), it follows that
$$
\mathrm{det} \Big(\frac{\partial \psi}{\partial s}, \frac{\partial \psi}{\partial t},\frac{\partial \psi}{\partial u}\Big)\Big |_{(0, 0, 0)} \neq 0.
$$
By the inverse function theorem,
it is obtained that there exist $0< a_{1}<a$, $0< b_{1}<b$, and $0< c_{1}<c$
with $a_{1}, b_{1}, c_{1}< \min \{\hat{t}_{1}, \hat{t}_{2}, \hat{t}_{3}\}$,
such that the function $\psi(s,t, u)$ is a diffeomorphism between
$V=(-a_{1},a_{1})\times (-b_{1},b_{1})\times (-c_{1},c_{1})$ and $U=\psi(V)$.
Hence, for each point $(u_{1},u_{2},u_{3})\in U$,
a point $(\overline{s}, \overline{t},\overline{u})\in V$ exists such that
$$(u_{1},u_{2},u_{3})=\psi(\overline{s}, \overline{t},\overline{u})\in \Gamma,$$
which, together with the assertion (a), implies that
$U\subset \Gamma\subset \Omega(z_{0}, \omega)$\ a.s.
Note that $U$ is a open subset of $\mathcal{K}\subset \mathbb{R}^{3}_{+}$.
Therefore, a almost surely finite stopping time $\widetilde{T}$
with respect to filtration $\mathcal{F}_{t}$
exists such that $(S(\widetilde{T}),I(\widetilde{T}),R(\widetilde{T}))\in U$\ a.s.
And because $\Gamma$ is a positive invariant set for the solutions of system (\ref{2.2})
and $U\subset \Gamma$, we can get that for almost all $\omega\in \widetilde{\Omega}$,
$(S(t),I(t),R(t))\in \Gamma$ is valid for all $t>\widetilde{T}$.
This also implies that $\Omega(z_{0}, \omega)\subset \overline{\Gamma}$\ a.s.
Thus, combining with the claim (a) of this theorem,
one can get $\overline{\Gamma}=\Omega(z_{0}, \omega) $\ a.s.
This completes the proof of Theorem \ref{theorem4.1.}.\ \ $\Box$

\end{pot1}

\subsection{Proof of Theorem \ref{theorem5.1.}}

To prove Theorem \ref{theorem5.1.}, we need the following two lemmas.

\begin{lem}\label{lemma5.1.}
For any $\varepsilon>0$, there exists a $\overline{T}=\overline{T}(\varepsilon)>0$
and a subset $A \in \mathcal{F}_{\infty}$ with $\mathbb{P}(A)>1-\varepsilon$
such that for any $t>\overline{T}$ and $S_{0}\in [0, \Lambda/\mu]$,
\begin{eqnarray}\label{5.1.1}
  \Bigg| \frac{1}{t}\int_{0}^{t} \Big(\beta_{r(s, \omega)}\widetilde{S}(s, S_{0})G'(0)
  -(\mu+\alpha+\delta)\Big)ds
  - \sum_{e\in \mathcal{M}}\pi_{e}B(e) \Bigg| &<& \varepsilon
\end{eqnarray}
holds for every $\omega\in A$, where $\widetilde{S}(s, S_{0})$
is the solution to system
\begin{eqnarray}\label{defi}
  \frac{d \widetilde{S}(t)}{dt} &=& \Lambda-\mu \widetilde{S}(t)
\end{eqnarray}
with the initial condition $\widetilde{S}(0)=S_{0}$.
\end{lem}

\begin{pf}
Consider the equation
\begin{eqnarray*}
  \frac{d \widetilde{S}(t)}{dt} &=& \Lambda-\mu \widetilde{S}(t).
\end{eqnarray*}
Let $\varepsilon_{1}$ be any positive constant satisfying
$\varepsilon_{1}< \varepsilon/\big(8G'(0)\sum_{e\in \mathcal{M}}\pi_{e}\beta_{e}\big)$.
By the compactness of the interval $[0, \Lambda/\mu]$
and the continuous dependence of solutions on initial conditions,
there exists a positive number $T_{1}=T_{1}(\varepsilon_{1})=T_{1}(\varepsilon)$
such that
\begin{eqnarray}\label{5.1.2}
  \Big| \widetilde{S}(t, S_{0})-\frac{\Lambda}{\mu} \Big| &<& \varepsilon_{1}
\end{eqnarray}
for all $t> T_{1}$ and $S_{0}\in [0, \Lambda/\mu]$.
Let us first fix a $S_{0}^{*}\in [0, \Lambda/\mu]$.
For every $\omega \in \widetilde{\Omega}$ and $t> T_{1}$,
we then have
\begin{eqnarray}\label{5.1.3}
  && \Bigg| \frac{1}{t}\int_{0}^{t} \Big(\beta_{r(s, \omega)}\widetilde{S}(s, S_{0}^{*})G'(0)
  -(\mu+\alpha+\delta)\Big)ds
  - \sum_{e\in \mathcal{M}}\pi_{e}B(e) \Bigg| \nonumber \\ [+6pt]
  &\leq&  \Bigg| \frac{1}{t}\int_{0}^{t} \beta_{r(s, \omega)}\Big(\widetilde{S}(s, S_{0}^{*})-\frac{\Lambda}{\mu}\Big)G'(0)ds \Bigg| \nonumber \\[+6pt]
  && +  \Bigg| \frac{1}{t}\int_{0}^{t} \Big(\beta_{r(s, \omega)}\frac{\Lambda G'(0)}{\mu}
  -(\mu+\alpha+\delta)\Big)ds - \sum_{e\in \mathcal{M}}\pi_{e}B(e) \Bigg| \nonumber \\[+6pt]
   &\leq&  \frac{1}{t}\int_{0}^{T_{1}} \beta_{r(s, \omega)}\Bigg| \widetilde{S}(s, S_{0}^{*})-\frac{\Lambda}{\mu}\Bigg|G'(0)ds
   + \frac{1}{t}\int_{T_{1}}^{t} \beta_{r(s, \omega)}\Bigg| \widetilde{S}(s, S_{0}^{*})-\frac{\Lambda}{\mu}\Bigg|G'(0)ds \nonumber \\[+6pt]
   && +  \Bigg| \frac{1}{t}\int_{0}^{t} \Big(\beta_{r(s, \omega)}\frac{\Lambda G'(0)}{\mu}
  -(\mu+\alpha+\delta)\Big)ds - \sum_{e\in \mathcal{M}}\pi_{e}B(e) \Bigg|.
\end{eqnarray}
Since $\beta_{r(s, \omega)}\Big| \widetilde{S}(s, S_{0}^{*})-\frac{\Lambda}{\mu}\Big|G'(0)
\leq \frac{\Lambda \beta^{M}G'(0) }{\mu}$ for all $s\in [0, T_{1}]$ and $\omega \in \widetilde{\Omega}$,
there exists a $T_{2}=T_{2}(\varepsilon)\geq T_{1}$ such that for $t> T_{2}$,
\begin{eqnarray}\label{5.1.4}
\frac{1}{t}\int_{0}^{T_{1}} \beta_{r(s, \omega)}\Bigg| \widetilde{S}(s, S_{0}^{*})
-\frac{\Lambda}{\mu}\Bigg|G'(0)ds  &<& \frac{\varepsilon}{3}.
\end{eqnarray}
Moreover, it is obtained from (\ref{5.1.2}) that
\begin{eqnarray*}
\frac{1}{t}\int_{T_{1}}^{t} \beta_{r(s, \omega)}\Bigg| \widetilde{S}(s, S_{0}^{*})
-\frac{\Lambda}{\mu}\Bigg|G'(0)ds  &\leq& \frac{\varepsilon_{1}G'(0)}{t-T_{1}}\int_{T_{1}}^{t} \beta_{r(s, \omega)}ds.
\end{eqnarray*}
From the Birkhoff Ergodic theorem,
there exists a subset $A_{1}\in \mathcal{F}_{\infty}$
satisfying $\mathbb{P}(A_{1})=1$ such that
for each $\omega\in A_{1}$, a $T_{3}=T_{3}(\omega, \varepsilon)\geq T_{1}$
exists, and when $t> T_{3}$,
\begin{eqnarray}\label{5.1.5}
\frac{1}{t}\int_{T_{1}}^{t} \beta_{r(s, \omega)}\Bigg| \widetilde{S}(s, S_{0}^{*})
-\frac{\Lambda}{\mu}\Bigg|G'(0)ds  &\leq&
2\varepsilon_{1}G'(0)\sum_{e\in \mathcal{M}}\pi_{e}\beta_{e}<\frac{\varepsilon}{3}
\end{eqnarray}
and
\begin{eqnarray}\label{5.1.6}
\Bigg| \frac{1}{t}\int_{0}^{t} \Big(\beta_{r(s, \omega)}\frac{\Lambda G'(0)}{\mu}
  -(\mu+\alpha+\delta)\Big)ds - \sum_{e\in \mathcal{M}}\pi_{e}B(e) \Bigg|
  &<&  \frac{\varepsilon}{3}.
\end{eqnarray}
Hence, substituting (\ref{5.1.4}), (\ref{5.1.5}) and (\ref{5.1.6}) into (\ref{5.1.3}), it is obtained that
for $t> T_{4}(\omega, \varepsilon)= \max \{T_{2}, T_{3}\}$,
\begin{eqnarray*}
 \Bigg| \frac{1}{t}\int_{0}^{t} \Big(\beta_{r(s, \omega)}\widetilde{S}(s, S_{0}^{*})G'(0)
  -(\mu+\alpha+\delta)\Big)ds
  - \sum_{e\in \mathcal{M}}\pi_{e}B(e) \Bigg| &<& \varepsilon
\end{eqnarray*}
holds for each $\omega \in A_{1}$. This implies
\begin{eqnarray*}
 \lim_{t\rightarrow +\infty}\frac{1}{t}\int_{0}^{t} \Big(\beta_{r(s)}\widetilde{S}(s, S_{0}^{*})G'(0)
  -(\mu+\alpha+\delta)\Big)ds&=& \sum_{e\in \mathcal{M}}\pi_{e}B(e)\ \ \ \ \ \ \ a.s.
\end{eqnarray*}
By the Egoroff theorem, there exists a $T_{5}=T_{5}(\varepsilon)>0$
and a subset $A_{2}\in \mathcal{F}_{\infty}$ with $\mathbb{P}(A_{2})>1-\varepsilon /2$
such that for any $t>T_{5}$ and $\omega \in A_{2}$,
\begin{eqnarray*}
 \Bigg| \frac{1}{t}\int_{0}^{t} \Big(\beta_{r(s, \omega)}\widetilde{S}(s, S_{0}^{*})G'(0)
  -(\mu+\alpha+\delta)\Big)ds
  - \sum_{e\in \mathcal{M}}\pi_{e}B(e) \Bigg| &<& \frac{\varepsilon}{2}.
\end{eqnarray*}
From (\ref{5.1.2}), it can be seen that for any $S_{1}, S_{2}\in [0, \Lambda/\mu]$
\begin{eqnarray*}
  \Big| \widetilde{S}(t, S_{1})-\widetilde{S}(t, S_{2}) \Big| &<& 2 \varepsilon_{1}
  \ \ \ \ \mathrm{for}\ \mathrm{all}\ t> T_{1},
\end{eqnarray*}
which, together with the preceding similar arguments, implies that
there exists a $T_{6}=T_{6}(\varepsilon)> T_{1}$
and a subset $A_{3}\in \mathcal{F}_{\infty}$ with
$\mathbb{P}(A_{3})> 1- \varepsilon /2$ satisfying that
\begin{eqnarray*}
\frac{1}{t}\int_{0}^{t} \beta_{r(s, \omega)}
 \Bigg| \widetilde{S}(s, S_{0})-\widetilde{S}(s, S_{0}^{*})\Bigg|G'(0)ds
&<& \frac{\varepsilon}{2}
\end{eqnarray*}
for any $S_{0}\in [0, \Lambda/\mu]$, $t> T_{6}$ and $\omega\in A_{3}$.
Let $A=A_{2}\bigcap A_{3}$, then $\mathbb{P}(A)>1-\varepsilon$.
For any $t> \overline{T}(\varepsilon)=\max \{T_{5}, T_{6}\}$
and $S_{0}\in [0, \Lambda/\mu]$,
\begin{eqnarray*}
  && \Bigg| \frac{1}{t}\int_{0}^{t} \Big(\beta_{r(s, \omega)}\widetilde{S}(s, S_{0})G'(0)
  -(\mu+\alpha+\delta)\Big)ds
  - \sum_{e\in \mathcal{M}}\pi_{e}B(e) \Bigg| \\ [+6pt]
  &\leq& \frac{1}{t}\int_{0}^{t} \beta_{r(s, \omega)}
 \Bigg| \widetilde{S}(s, S_{0})-\widetilde{S}(s, S_{0}^{*})\Bigg|G'(0)ds\\[+6pt]
 && +\Bigg| \frac{1}{t}\int_{0}^{t} \Big(\beta_{r(s, \omega)}\widetilde{S}(s, S_{0}^{*})G'(0)
  -(\mu+\alpha+\delta)\Big)ds
  - \sum_{e\in \mathcal{M}}\pi_{e}B(e) \Bigg|\\[+6pt]
  &\leq& \frac{\varepsilon}{2}+\frac{\varepsilon}{2}<\varepsilon
\end{eqnarray*}
holds for each $\omega \in A$.
This completes the proof of Lemma \ref{lemma5.1.}.\ \ $\Box$

\end{pf}

\begin{lem}\label{lemma5.2.}
For any $\varepsilon>0$, there exists a
$\gamma=\gamma(\varepsilon)\in (0, \Lambda/(\mu+\alpha)]$
such that the following statements are valid:

(i) there exist two positive constants $\widetilde{T}_{1}=\widetilde{T}_{1}(\varepsilon)$
and $\widetilde{T}_{2}=\widetilde{T}_{2}(\varepsilon)$,
such that
for any $t> \widetilde{T}_{1}+\widetilde{T}_{2}$,
if $I(s,\omega, z_{0})< \gamma$ for all $s\in [0,t]$,
then
\begin{eqnarray*}
  \Big| S(s,\omega, z_{0})-\widetilde{S}(s, S_{0}) \Big|<\frac{\varepsilon}{3L}
\end{eqnarray*}
holds for all $s\in [\widetilde{T}_{1}+\widetilde{T}_{2},t]$;

(ii) for any $t>0$, if $I(s,\omega, z_{0})< \gamma$ for all $s\in [0,t]$,
then
\begin{eqnarray*}
\Bigg|  \frac{G(I(s,\omega, z_{0}))}{I(s,\omega, z_{0})}-G'(0)\Bigg|<\frac{\varepsilon}{3L}
\end{eqnarray*}
holds for all $s\in [0,t]$,
where $L=\max\{\beta^{M}G'(0), (\Lambda\beta^{M})/\mu\}$ and
$\widetilde{S}(s, S_{0})$
is the solution of system (\ref{defi})
with the initial condition $\widetilde{S}(0)=S_{0}\in [0, \Lambda/\mu]$.
\end{lem}

\begin{pf}
(i) Let $\varepsilon_{1}=\mu\varepsilon/6L$.
Consider the following equation
\begin{eqnarray}\label{uu1}
\frac{du^{1}(t)}{dt} &=& -\mu u^{1}(t)+\varepsilon_{1}
\end{eqnarray}
with the initial value $u^{1}(0)=u^{1}_{0}\in \mathbb{R}$.
It is easy to see that the system (\ref{uu1})
has a globally asymptotically stable equilibrium $u_{1}^{*}=\varepsilon_{1} / \mu$.
By Lemma \ref{lemma-global}, it follows that
a constant $\widetilde{T}_{1}>0$ exists
such that $u^{1}(t)\leq 2 \varepsilon_{1} / \mu$ for all $t\geq \widetilde{T}_{1}$ and
$u^{1}_{0}\in [-\Lambda/\mu, \Lambda/\mu]$.

Suppose that for a fixed constant $\gamma_{1}>0$,
$I(s,\omega, z_{0})< \gamma_{1}$ for all $s\in [0,t]$
with $t> \widetilde{T}_{1}+\widetilde{T}_{2}$,
where $\gamma_{1}=\gamma_{1}(\varepsilon)$ and $\widetilde{T}_{2}=\widetilde{T}_{2}(\gamma_{1})$
will be determined later.
It is obtained from the third equation of system (\ref{2.2}) that
\begin{eqnarray*}
\frac{d R(t,\omega, z_{0})}{dt} &=& \delta I(t,\omega, z_{0})-(\mu+\lambda)R(t,\omega, z_{0})\\[+6pt]
&\leq& \delta \gamma_{1}-(\mu+\lambda)R(t,\omega, z_{0})
\end{eqnarray*}
for all $s\in [0,t]$.
Moreover, for the system
\begin{eqnarray}\label{uu2}
\frac{du^{2}(t)}{dt} &=& \delta \gamma_{1}-(\mu+\lambda)u^{2}(t)
\end{eqnarray}
with the initial value $u^{2}(0)=\Lambda/\mu$,
one can easily find a constant $\widetilde{T}_{2}=\widetilde{T}_{2}(\gamma_{1})>0$
such that $u^{2}(t)<2\delta \gamma_{1}/ (\mu+\lambda)$ if $t\geq \widetilde{T}_{2}$.
By the comparison theorem, it follows from
(\ref{uu2}) that $R(t,\omega, z_{0})<2\delta \gamma_{1}/ (\mu+\lambda)$
when $t\geq \widetilde{T}_{2}$.
Define the functions $F_{e}(S,I,R)=\Lambda-\mu S+\lambda R-\beta_{e}S G(I)$,
$e\in \mathcal{M}$.
Due to the continuity of $F_{e}$ with respect to variables $I$ and $R$,
one can find a sufficiently small $\gamma_{1}=\gamma_{1}(\varepsilon)\in (0, \Lambda/(\mu+\alpha)]$
such that whenever $0\leq I(s,\omega, z_{0})\leq \gamma_{1} $
and $0\leq S(s,\omega, z_{0})\leq\Lambda/\mu$ for all $s\in [0,t]$,
$$\Lambda-\mu S(s,\omega, z_{0})-\varepsilon_{1}\leq \frac{d S(s,\omega, z_{0})}{ds}
\leq \Lambda-\mu S(s,\omega, z_{0})+\varepsilon_{1}$$
for all $s\in [\widetilde{T}_{2},t]$.
Introduce the following equation
\begin{eqnarray}\label{defi2}
  \frac{du(s)}{ds} &=& \Lambda-\mu u(s)-\varepsilon_{1},
\end{eqnarray}
with the initial condition $u(\widetilde{T}_{2})=S(\widetilde{T}_{2},\omega, z_{0})$.
From the comparison theorem,
it is easy to see that
\begin{eqnarray}\label{5.1.7}
  S(s,\omega, z_{0}) &\geq& u(s,u(\widetilde{T}_{2}))\ \ \ \mathrm{for} \ \mathrm{all}\ \ s\in [\widetilde{T}_{2},t].
\end{eqnarray}
From (\ref{defi}) and (\ref{defi2}), it is obtained that for $s\in [\widetilde{T}_{2},t]$
\begin{eqnarray*}
  \frac{d}{ds}\big(\widetilde{S}(s)-u(s)\big)&=&-\mu \big(\widetilde{S}(s)-u(s)\big)+\varepsilon_{1}
\end{eqnarray*}
with the initial value $\widetilde{S}(\widetilde{T}_{2})-u(\widetilde{T}_{2})
\in [-\Lambda/\mu, \Lambda/\mu]$.
Hence,
it is obtained from (\ref{uu1}) and (\ref{5.1.7}) that for $s\in [\widetilde{T}_{1}+\widetilde{T}_{2},t]$,
\begin{eqnarray*}
 \widetilde{S}(s, S_{0})-S(s, \omega, z_{0})&=& \big( \widetilde{S}(s, S_{0})
 -u(s, u(\widetilde{T}_{2}))\big)+ \big(u(s, u(\widetilde{T}_{2}))-S(s,\omega, z_{0}) \big)\\[+6pt]
   &\leq& \frac{2\varepsilon_{1}}{\mu}< \frac{\varepsilon}{3L}.
\end{eqnarray*}
Using the similar arguments, we also have
\begin{eqnarray*}
S(s, \omega, z_{0})-  \widetilde{S}(s, S_{0})&<& \frac{\varepsilon}{3L}
\ \ \ \mathrm{for} \ \mathrm{all}\ \ s\in [\widetilde{T}_{1}+\widetilde{T}_{2},t].
\end{eqnarray*}
Hence, for all $s\in [\widetilde{T}_{1}+\widetilde{T}_{2},t]$, we have
\begin{eqnarray*}
\Big |S(s, \omega, z_{0})-  \widetilde{S}(s, S_{0})\Big |
&<& \frac{\varepsilon}{3L}.
\end{eqnarray*}

(ii) On the other hand, by the continuity of the function $g(I)=G(I)/I$
with respect to variable $I$
and $\lim_{I\rightarrow 0^{+}}g(I)=G'(0)$,
it follows that there exists a sufficiently small $\gamma_{2}=\gamma_{2}(\varepsilon)\in (0, \Lambda/(\mu+\alpha)]$,
such that whenever $0\leq I(s,\omega, z_{0})\leq \gamma_{2}$,
\begin{eqnarray*}
\Bigg | \frac{G(I(s,\omega, z_{0}))}{I(s,\omega, z_{0})}-G'(0)\Bigg |&<& \frac{\varepsilon}{3L}.
\end{eqnarray*}

Let us choose $\gamma=\gamma(\varepsilon)=\min \{\gamma_{1},\gamma_{2}\}$.
This completes the proof of Lemma \ref{lemma5.2.}.\ \ $\Box$

\end{pf}

\begin{pot}
Let $\varepsilon$ be any positive constant with
$\varepsilon< (\sum_{e\in \mathcal{M}}\pi_{e}B(e))/2$.
Consider the process $((S(t), I(t), R(t)), r(t))$
on a larger state space
\begin{eqnarray*}
   \overline{X}= \widetilde{\mathcal{K}}\times \mathcal{M},
\end{eqnarray*}
where $\widetilde{\mathcal{K}}=\overline{\mathcal{K}}\backslash \{(S,I,R)\in \mathbb{R}^{3}_{+}: I=0\}$ and $\overline{\mathcal{K}}$ is the closure of $\mathcal{K}$.
Note that the process $((S(t), I(t), R(t)), r(t))$ is a time-homogeneous Markov process with
the Feller property. According to Theorem \ref{invariant-exist},
one can get the existence of an invariant
probability measure $\nu^{*}$ for the process
$((S(t), I(t), R(t)), r(t))$ on the state space $\overline{X}$,
provided that a compact set $O \subset \overline{X}$ exists
such that
\begin{eqnarray}\label{dis-existence}
\liminf_{t\rightarrow +\infty} \frac{1}{t}\int_{0}^{t}
\Bigg( \int_{\overline{X}}\mathbb{P}(s, x, O)\nu(dx)\Bigg)ds &=&
\liminf_{t\rightarrow +\infty} \frac{1}{t}\int_{0}^{t}
\mathbb{P}(s, x_{0}, O)ds> 0,
\end{eqnarray}
for some initial distribution $\nu=\delta_{x_{0}}$ with $x_{0}\in \overline{X}$, where
$\mathbb{P}(s, \cdot, \cdot)$ is the transition probability function
and $\delta_{\cdot}$ is the Dirac function.
On the other hand, it is easy to see that
for any initial value $x_{0}=\big(z_{0}, r(0)\big)\in X$,
when $t>0$ the solution $(S(t), I(t), R(t))$
of system (\ref{2.2}) do not reach the boundary $\partial \widetilde{\mathcal{K}}$ of the region
$\widetilde{\mathcal{K}}$ under the condition that $\mathcal{R}_{0}>1$.
Hence, we have $\nu^{*}(\partial \widetilde{\mathcal{K}} \times \mathcal{M})=0$,
which obviously implies that
$\nu^{*}$ is also the invariant
probability measure of the process
$((S(t), I(t), R(t)), r(t))$ on the state space $X$.
Consequently, to complete the proof of Theorem \ref{theorem5.1.},
it is sufficient to find a compact set $O \subset \overline{X}$
satisfying (\ref{dis-existence}).

Fix a $T> \max \{\overline{T}, \widetilde{T}_{1}+\widetilde{T}_{2},
6\Lambda L(\widetilde{T}_{1}+\widetilde{T}_{2})/\mu \varepsilon\}$, and set $\chi_{n}(\omega)=\mathbf{1}_{A}(\theta^{nT}\omega)$,
$n\in \mathbb{N}_{0}$, where $A$ is as in Lemma \ref{lemma5.1.},
$\mathbf{1}_{A}(\cdot)$ is the indicator function, and the shift operators $\theta^{t}$
on $\widetilde{\Omega}$ are defined by
$$(\theta^{t}\omega)_{s}=\omega_{s+t},\ \ \ s,t\geq 0.$$
Let
\begin{eqnarray*}
\chi_{n}^{1}(\omega)=
\left\{
  \begin{array}{l}
 1,\ \  \  \mathrm{if} \ \chi_{n}(\omega)=1 \ \mathrm{and} \ I(t,\omega,z_{0})
 <\gamma \\
 \ \ \ \ \ \ \mathrm{for} \ \mathrm{any} \ t\in[nT, (n+1)T],    \\ [+6pt]
 0,\ \ \ \mathrm{otherwise}
\end{array} \right.
\end{eqnarray*}
and
\begin{eqnarray*}
\chi_{n}^{2}(\omega)=
\left\{
  \begin{array}{l}
 1,\ \ \  \mathrm{if} \ \chi_{n}(\omega)=1 \ \mathrm{and} \ \mathrm{there} \
 \mathrm{exists}\  \mathrm{some}
 \ t\in[nT, (n+1)T]  \\
 \ \ \ \ \ \   \mathrm{such} \  \mathrm{that} \ I(t,\omega,z_{0})\geq \gamma,    \\ [+6pt]
 0,\ \ \ \mathrm{otherwise},
\end{array} \right.
\end{eqnarray*}
where $n\in \mathbb{N}_{0}$. It is easy to see that $\chi_{n}(\omega)=\chi_{n}^{1}(\omega)
+\chi_{n}^{2}(\omega)$. Denote $\chi_{n}^{3}(\omega)=1-\chi_{n}(\omega)$.
Since $r(t+s, \omega)=r(t, \theta^{s}\omega)$ for any $s,t>0$, we have
\begin{eqnarray*}
z(t+s, \omega, z_{0}) &=& z(t, \theta^{s}\omega, z(s,\omega,z_{0}))
\ \ \ \mathrm{for}\ \mathrm{any}\ s,t>0.
\end{eqnarray*}
This implies that if $\chi_{n}^{1}(\omega)=1$ then $I(t+nT,\omega,z_{0})< \gamma$
holds for all $0\leq t\leq T$,
i.e., $I(t,\theta^{nT}\omega,z(nT,\omega,z_{0}))< \gamma$
holds for all $0\leq t\leq T$.
In addition, $\chi_{n}^{1}(\omega)=1$ implies $\chi_{n}(\omega)=1$,
i.e., $\theta^{nT}\omega \in A$. Hence, it is obtained from
Lemmas \ref{lemma5.1.} and \ref{lemma5.2.} that
\begin{eqnarray}\label{5.1.11}
\Bigg| \frac{1}{T}\int_{0}^{T} \Big(\beta_{r(t, \theta^{nT}\omega)}
\widetilde{S}\big(t, S(nT,\omega,z_{0})\big)G'(0)
  -(\mu+\alpha+\delta)\Big)dt
  - \sum_{e\in \mathcal{M}}\pi_{e}B(e) \Bigg| &<& \varepsilon,
\end{eqnarray}
\begin{eqnarray}\label{5.1.12-1}
  \frac{1}{T-(\widetilde{T}_{1}+\widetilde{T}_{2})}
  \int_{\widetilde{T}_{1}+\widetilde{T}_{2}}^{T} \Bigg| S(t,\theta^{nT}\omega, z(nT,\omega,z_{0}))
  -\widetilde{S}(t, S(nT,\omega,z_{0})) \Bigg|dt &<&\frac{\varepsilon}{3L}
\end{eqnarray}
and
\begin{eqnarray}\label{5.1.12-2}
  \frac{1}{T}\int_{0}^{T}\Bigg|  \frac{G\big(I(t,\theta^{nT}\omega, z(nT,\omega,z_{0}))\big)}{I(t,\theta^{nT}\omega, z(nT,\omega,z_{0}))}-G'(0)\Bigg|dt &<&\frac{\varepsilon}{3L}.
\end{eqnarray}
For the second equation of system (\ref{2.2}),
we compute that
\begin{eqnarray*}
&& \frac{1}{T}\Big(\ln I((n+1)T, \omega, z_{0})-\ln I(nT, \omega, z_{0}) \Big) \\ [+6pt]
&=& \frac{1}{T}\int_{nT}^{(n+1)T}\Bigg( \beta_{r(t,\omega)}\frac{S(t,\omega,z_{0})G(I(t,\omega,z_{0}))}
{I(t,\omega,z_{0})}-(\mu+\alpha+\delta)\Bigg)dt \\[+6pt]
&=& \frac{1}{T}\int_{0}^{T}\Bigg( \beta_{r(t+nT,\omega)}\frac{S(t+nT,\omega,z_{0})G(I(t+nT,\omega,z_{0}))}
{I(t+nT,\omega,z_{0})}-(\mu+\alpha+\delta)\Bigg)dt \\[+6pt]
&=& \frac{1}{T}\int_{0}^{T}\Bigg( \beta_{r(t,\theta^{nT}\omega)}\frac{S(t,\theta^{nT}\omega,z(nT,\omega,z_{0}))
G(I(t,\theta^{nT}\omega,z(nT,\omega,z_{0})))}
{I(t,\theta^{nT}\omega,z(nT,\omega,z_{0}))}-(\mu+\alpha+\delta)\Bigg)dt \\[+6pt]
&=& \frac{1}{T}\int_{0}^{T} \Big(\beta_{r(t, \theta^{nT}\omega)}
\widetilde{S}\big(t, S(nT,\omega,z_{0})\big)G'(0)-(\mu+\alpha+\delta)\Big)dt
+\frac{1}{T}\int_{0}^{T} \beta_{r(t,\theta^{nT}\omega)}\\[+6pt]
&&\Big( S(t,\theta^{nT}\omega, z(nT,\omega,z_{0}))
  -\widetilde{S}(t, S(nT,\omega,z_{0})) \Big)
  \frac{G(I(t,\theta^{nT}\omega,z(nT,\omega,z_{0})))}{I(t,\theta^{nT}\omega,z(nT,\omega,z_{0}))}dt  \\[+6pt]
&& +\frac{1}{T}\int_{0}^{T} \beta_{r(t,\theta^{nT}\omega)}\widetilde{S}(t, S(nT,\omega,z_{0}))
  \Bigg(\frac{G(I(t,\theta^{nT}\omega,z(nT,\omega,z_{0})))}{I(t,\theta^{nT}\omega,z(nT,\omega,z_{0}))}-G'(0) \Bigg)dt  \\[+6pt]
&\geq& \frac{1}{T}\int_{0}^{T} \Big(\beta_{r(t, \theta^{nT}\omega)}
\widetilde{S}\big(t, S(nT,\omega,z_{0})\big)G'(0)-(\mu+\alpha+\delta)\Big)dt\\[+6pt]
&&\ \ \ \ \ \ -\frac{L}{T}\int_{0}^{\widetilde{T}_{1}+\widetilde{T}_{2}}\Bigg( \Bigg| S(t,\theta^{nT}\omega, z(nT,\omega,z_{0}))
  -\widetilde{S}(t, S(nT,\omega,z_{0})) \Bigg|\Bigg)dt  \\[+6pt]
  &&\ \ \ \ \ \ -\frac{L}{T-(\widetilde{T}_{1}+\widetilde{T}_{2})}
  \int_{\widetilde{T}_{1}+\widetilde{T}_{2}}^{T}\Bigg( \Bigg| S(t,\theta^{nT}\omega, z(nT,\omega,z_{0}))
  -\widetilde{S}(t, S(nT,\omega,z_{0})) \Bigg|\Bigg)dt  \\[+8pt]
  &&\ \ \ \ \ \ -\frac{L}{T}\int_{0}^{T}\Bigg(\Bigg|  \frac{G\big(I(t,\theta^{nT}\omega, z(nT,\omega,z_{0}))\big)}{I(t,\theta^{nT}\omega, z(nT,\omega,z_{0}))}-G'(0)\Bigg|\Bigg)dt.
\end{eqnarray*}
Substituting (\ref{5.1.11}), (\ref{5.1.12-1}) and (\ref{5.1.12-2}) into the right side of the above inequality,
it is obtained that when $\chi_{n}^{1}(\omega)=1$, we have
\begin{eqnarray}\label{5.1.13}
\frac{1}{T}\Big(\ln I((n+1)T, \omega, z_{0})-\ln I(nT, \omega, z_{0}) \Big)
\geq  \sum_{e\in \mathcal{M}}\pi_{e}B(e)-2\varepsilon.
\end{eqnarray}
Let
$$\Theta=\max \Bigg\{ \Bigg|\beta_{e}\frac{SG(I)}{I}-(\mu+\alpha+\delta)  \Bigg|:
 0\leq S,I \leq \frac{\Lambda}{\mu}, e\in \mathcal{M}\Bigg\}.$$
If $\chi_{n}^{2}(\omega)=1$ or $\chi_{n}^{3}(\omega)=1$, then
\begin{eqnarray}\label{5.1.14}
&&\frac{1}{T}\Big(\ln I((n+1)T, \omega, z_{0})-\ln I(nT, \omega, z_{0}) \Big) \nonumber \\[+6pt]
&=& \frac{1}{T}\int_{nT}^{(n+1)T}\Bigg( \beta_{r(t,\omega)}\frac{S(t,\omega,z_{0})G(I(t,\omega,z_{0}))}{I(t,\omega,z_{0})}
-(\mu+\alpha+\delta) \Bigg)dt \nonumber \\[+6pt]
 &\geq& -\Theta.
\end{eqnarray}
Summing up both sides of (\ref{5.1.13}) and (\ref{5.1.14}) respectively yields
$$\frac{1}{T}\Big(\ln I((n+1)T, \omega, z_{0})-\ln I(nT, \omega, z_{0}) \Big)
\geq \Big(\sum_{e\in \mathcal{M}}\pi_{e}B(e)-2\varepsilon \Big)
\chi_{n}^{1}(\omega)-\Theta \Big(\chi_{n}^{2}(\omega)+\chi_{n}^{3}(\omega)\Big),$$
which implies
\begin{eqnarray*}
 && \frac{1}{kT}\Big(\ln I(kT, \omega, z_{0})-\ln I(0) \Big) \\[+6pt]
   &&\ \ \ \ \ \ \ \ \ \ \ \ \ \ \ \ \geq \frac{1}{k}\Bigg(\Big(\sum_{e\in \mathcal{M}}\pi_{e}B(e)-2\varepsilon \Big)
\sum_{n=0}^{k-1}\chi_{n}^{1}(\omega)-\Theta \sum_{n=0}^{k-1}\Big(\chi_{n}^{2}(\omega)+\chi_{n}^{3}(\omega)\Big)\Bigg)
\end{eqnarray*}
holds for any $k\in \mathbb{N}$. Noting that $I(t)\leq \Lambda/\mu$ for all
$t\geq0$, we then have
\begin{eqnarray}\label{5.1.15}
\limsup_{k\rightarrow +\infty}\frac{1}{k}\Bigg(\Big(\sum_{e\in \mathcal{M}}\pi_{e}B(e)-2\varepsilon \Big)
\sum_{n=0}^{k-1}\chi_{n}^{1}(\omega)-\Theta \sum_{n=0}^{k-1}\Big(\chi_{n}^{2}(\omega)+\chi_{n}^{3}(\omega)\Big)\Bigg) &\leq& 0.
\end{eqnarray}
Using the arguments similar to that given in Theorem 3.1 of \cite{Yingang2014},
we have
$$\lim_{k \rightarrow +\infty}\frac{1}{k}\sum_{n=0}^{k-1}\chi_{n}(\omega)=\mathbb{P}(A)\ \ \ \ a.s.,$$
which implies
\begin{eqnarray}\label{5.1.16}
\lim_{k \rightarrow +\infty}\frac{1}{k}\sum_{n=0}^{k-1}\chi_{n}^{3}(\omega) &=& 1- \mathbb{P}(A)\leq
\varepsilon\ \ \ \ a.s.,
\end{eqnarray}
and
\begin{eqnarray}\label{5.1.17}
\lim_{k \rightarrow +\infty}\frac{1}{k}\sum_{n=0}^{k-1}\Big(\chi_{n}^{1}(\omega)+\chi_{n}^{2}(\omega) \Big)
&=& \lim_{k \rightarrow +\infty}\frac{1}{k}\sum_{n=0}^{k-1}\chi_{n}(\omega)\geq 1- \varepsilon \ \ \ \ a.s.
\end{eqnarray}
Multiplying both sides of (\ref{5.1.16}) by $\Theta$ and multiplying both sides
of (\ref{5.1.17}) by $-\big( \sum_{e\in \mathcal{M}}\pi_{e}B(e)$ $-$ $2\varepsilon \big)$,
then adding to (\ref{5.1.15}) yields
\begin{eqnarray*}
&& \limsup_{k\rightarrow +\infty}\frac{1}{k}\Bigg[\Delta
\sum_{n=0}^{k-1}\chi_{n}^{1}(\omega)-\Theta \sum_{n=0}^{k-1}\chi_{n}^{2}(\omega)- \Delta\sum_{n=0}^{k-1}\Big(\chi_{n}^{1}(\omega)+\chi_{n}^{2}(\omega)\Big)\Bigg]\\[+6pt]
&=&  \limsup_{k\rightarrow +\infty}\frac{1}{k}\Big[-\big(\Theta+\Delta \big) \Big]\sum_{n=0}^{k-1}\chi_{n}^{2}(\omega)\\[+6pt]
&\leq& -\Delta(1-\varepsilon)+\Theta\varepsilon\ \ \ \ a.s.,
\end{eqnarray*}
where $\Delta=\sum_{e\in \mathcal{M}}\pi_{e}B(e)-2\varepsilon >0$.
Hence, it is obtained that for $\varepsilon$ sufficiently small,
\begin{eqnarray}\label{5.1.18}
\liminf_{k\rightarrow +\infty} \frac{1}{k}\sum_{n=0}^{k-1}\chi_{n}^{2}(\omega)
&\geq& \frac{\Delta(1-\varepsilon)-\Theta \varepsilon}{\Theta+\Delta}:=\kappa>0 \ \ \ \ a.s.
\end{eqnarray}
Denote by $D(s)$ the closure of the region
$$\Bigg\{(S,I,R)\in \mathbb{R}^{3}_{+}: \frac{\Lambda}{\mu+\alpha}\leq S+I+R\leq\frac{\Lambda}{\mu},
I \geq s, 0<s<\frac{\Lambda}{\mu}  \Bigg \}.$$
For the system (\ref{2.2}), one can find a $\overline{\gamma}=\overline{\gamma}(\gamma,T)>0$
satisfying $z(t,\omega,z_{0})\in D(\overline{\gamma})$ for all $t\in [0, T]$ and $z_{0}\in \mathcal{K}$,
provided that an $s\in [0, T]$ exists such that $z(s,\omega,z_{0})\in D(\gamma)$.
Hence, for each $\omega \in \widetilde{\Omega}$ with $\chi_{n}^{2}(\omega)=1$
and any $z_{0}\in \mathcal{K}$, we have
$$\int_{nT}^{(n+1)T}\mathbf{1}_{\{z(t,\omega,z_{0})\in D(\overline{\gamma})\}}dt=T,$$
which, combined with (\ref{5.1.18}), implies
$$\liminf_{t\rightarrow +\infty}\frac{1}{t}\int_{0}^{t}\mathbf{1}_{\{z(s,\omega,z_{0})\in D(\overline{\gamma})\}}ds\geq \kappa>0 \ \ \ \ a.s.$$
By Fatou's lemma, it then follows that
\begin{eqnarray*}
\liminf_{t\rightarrow +\infty}\frac{1}{t}\int_{0}^{t}\mathbb{P}
\Big\{ z(s,\omega,z_{0})\in D(\overline{\gamma}) \Big\}ds &\geq& \kappa>0.
\end{eqnarray*}
Let $O= D(\overline{\gamma})\times \mathcal{M}\subset \overline{X}$,
we then have
\begin{eqnarray*}
\liminf_{t\rightarrow +\infty}\frac{1}{t}\int_{0}^{t}\mathbb{P}
(s,x_{0}, O )ds=\liminf_{t\rightarrow +\infty}\frac{1}{t}\int_{0}^{t}\mathbb{P}
\Big\{ z(s,\omega,z_{0})\in D(\overline{\gamma}) \Big\}ds &\geq& \kappa>0.
\end{eqnarray*}
This completes the proof of Theorem \ref{theorem5.1.}.\ \ $\Box$

\end{pot}

\begin{rmk}
The method used in the proof of Theorem \ref{theorem5.1.}
is similar to that given in Theorem 3.1 of \cite{Yingang2014}.
In the proceeding of the proof,
it is most important to find a compact set $O \subset \overline{X}$
satisfying (\ref{dis-existence}).
Indeed, according to the results of
Theorem \ref{theorem3.2.},
there exists a simpler method for getting such compact set $O \subset \overline{X}$
as follows. By Theorem \ref{theorem3.2.}, we have
\begin{eqnarray*}
  \liminf_{t\rightarrow +\infty} \frac{1}{t}\int_{0}^{t}I(s) ds \geq
         \frac{\mu^{2}}{\beta^{M}(\mu \vartheta+\beta^{M}(G'(0))^{2})\Lambda}
         \sum_{e\in \mathcal{M}}\pi_{e}B(e):=\kappa_{1}>0 \ \ \ a.s.
\end{eqnarray*}
if $\mathcal{R}_{0}>1$. Since
\begin{eqnarray*}
\frac{1}{t}\int_{0}^{t}I(s) ds &=& \frac{1}{t}\int_{0}^{t}I(s)\cdot
\mathbf{1}_{\{I(s)< \frac{\kappa_{1}}{2}\}} ds+
\frac{1}{t}\int_{0}^{t}I(s)\cdot
\mathbf{1}_{\{I(s)\geq \frac{\kappa_{1}}{2}\}} ds \\[+6pt]
 &\leq& \frac{\kappa_{1}}{2}+\frac{\Lambda}{\mu}\cdot
 \frac{1}{t}\int_{0}^{t}
\mathbf{1}_{\{I(s)\geq \frac{\kappa_{1}}{2}\}} ds,
\end{eqnarray*}
we have
\begin{eqnarray}\label{remark-exi}
\liminf_{t\rightarrow +\infty} \frac{1}{t}\int_{0}^{t}
\mathbf{1}_{\{I(s)\geq \frac{\kappa_{1}}{2}\}} ds
&\geq& \frac{\mu\kappa_{1}}{2\Lambda}\ \ \ a.s.
\end{eqnarray}
By Fatou's lemma, it then follows from (\ref{remark-exi}) that
\begin{eqnarray*}
\liminf_{t\rightarrow +\infty} \frac{1}{t}\int_{0}^{t}
\mathbb{P}\Big\{I(s)\geq \frac{\kappa_{1}}{2}\Big\}ds &=&
\liminf_{t\rightarrow +\infty} \frac{1}{t}\int_{0}^{t}\mathbb{E}
\Big[\mathbf{1}_{\{I(s)\geq \frac{\kappa_{1}}{2}\}}\Big]ds \\[+6pt]
&\geq& \mathbb{E}\Bigg[\liminf_{t\rightarrow +\infty} \frac{1}{t}\int_{0}^{t}
\mathbf{1}_{\{I(s)\geq \frac{\kappa_{1}}{2}\}} ds \Bigg] \\[+6pt]
&\geq& \frac{\mu\kappa_{1}}{2\Lambda}.
\end{eqnarray*}
Let $O= D(\frac{\kappa_{1}}{2}) \times \mathcal{M}\subset \overline{X}$.
It then follows that
\begin{eqnarray*}
\liminf_{t\rightarrow +\infty}\frac{1}{t}\int_{0}^{t}\mathbb{P}
(s,x_{0}, O )ds&=&\liminf_{t\rightarrow +\infty}\frac{1}{t}\int_{0}^{t}\mathbb{P}
\Big\{ z(s,\omega,z_{0})\in D \Big\}ds \\[+6pt]
&=&\liminf_{t\rightarrow +\infty} \frac{1}{t}\int_{0}^{t}
\mathbb{P}\Big\{I(s)\geq \frac{\kappa_{1}}{2}\Big\}ds \\[+6pt]
&\geq& \frac{\mu\kappa_{1}}{2\Lambda},
\end{eqnarray*}
which implies that (\ref{dis-existence}) holds.

\end{rmk}

\subsection{Proof of Theorem \ref{theorem5.2.}}

\begin{pf}
We shall complete the proof of this theorem by three steps.

{\bf Step 1.} We shall prove that the process $((S(t), I(t), R(t)), r(t))$
is positive Harris recurrent.
Since the point $\overline{z}_{*} \in \Gamma$ satisfies
the condition (\textbf{H}),
without loss of generality, we can assume by the similar arguments
in the proof of (b) in Theorem \ref{theorem4.1.} that the function
\begin{eqnarray}\label{5.2.1}
\psi(s,t, u) &=& \pi_{\hat{t}_{i}}^{\hat{p}_{i}}\circ \cdots \circ \pi_{\hat{t}_{4}}^{\hat{p}_{4}}\circ
\pi_{u+\hat{t}_{3}}^{\hat{p}_{3}} \circ \pi_{t+\hat{t}_{2}}^{\hat{p}_{2}}\circ \pi_{s+\hat{t}_{1}}^{\hat{p}_{1}}(\overline{z}_{*})
\end{eqnarray}
is defined and continuously differentiable in
$(-c,c)^{3}\subset \mathbb{R}^{3}$ with
\begin{eqnarray}\label{5.2.2}
\mathrm{det} \Big(\frac{\partial \psi}{\partial s}, \frac{\partial \psi}{\partial t},\frac{\partial \psi}{\partial u}\Big) &\neq& 0
\end{eqnarray}
for any $(s,t,u)\in (-c,c)^{3}$, where the states $\hat{p}_{k}\in \mathcal{M}$,
$\hat{t}_{k}>0$, $k=1,\ldots,i$, and $c$ is some positive constant.
From the definition of the set $\Gamma$, without loss of generality we can assume that
there exist $\bar{p}_{1}, \ldots, \bar{p}_{j}\in \mathcal{M}$
and $\bar{t}_{1}, \ldots, \bar{t}_{j}>0$
such that $\overline{z}_{*}=\pi_{\bar{t}_{j}}^{\bar{p}_{j}}\circ \cdots \circ \pi_{\bar{t}_{1}}^{\bar{p}_{1}}(E^{*}_{1})$.
Because $E^{*}_{1}$ is the equilibrium of system (\ref{2.2})
in the state 1, for any $w_{1}>0$ we then have $\pi_{w_{1}}^{1}(E^{*}_{1})=E^{*}_{1}$.
Hence, $\overline{z}_{*}$ can be rewritten as
\begin{eqnarray}\label{5.2.3}
\overline{z}_{*} &=& \pi_{\bar{t}_{j}}^{\bar{p}_{j}}\circ \cdots
\circ \pi_{\bar{t}_{1}}^{\bar{p}_{1}}\circ \pi_{w_{1}}^{1}(E^{*}_{1}).
\end{eqnarray}
Substituting (\ref{5.2.3}) into (\ref{5.2.1}), it is obtained that for any $w_{1}>0$,
$$\psi_{w_{1}}(s,t, u)=\psi(s,t, u)=\widetilde{F}(\tilde{t})\circ
\pi_{u+\hat{t}_{3}}^{\hat{p}_{3}} \circ \pi_{t+\hat{t}_{2}}^{\hat{p}_{2}}\circ \pi_{s+\hat{t}_{1}}^{\hat{p}_{1}}
\circ \overline{F}(\bar{t})\circ \pi_{w_{1}}^{1}(E^{*}_{1}),
$$
where $\tilde{t}=(\hat{t}_{4},\ldots,\hat{t}_{i})$, $\bar{t}=(\bar{t}_{1},\ldots,\bar{t}_{j})$ and
$\widetilde{F}(\tilde{t})=\pi_{\hat{t}_{i}}^{\hat{p}_{i}}\circ \cdots \circ \pi_{\hat{t}_{4}}^{\hat{p}_{4}}$,
$\overline{F}(\bar{t})= \pi_{\bar{t}_{j}}^{\bar{p}_{j}}\circ \cdots
\circ \pi_{\bar{t}_{1}}^{\bar{p}_{1}}$.
For each $z \in \mathbb{R}^{3}_{+}$,
we introduce the function as follows
$$\Psi_{z, w_{1}}(s,t, u)=\widetilde{F}(\tilde{t})\circ
\pi_{u+\hat{t}_{3}}^{\hat{p}_{3}} \circ \pi_{t+\hat{t}_{2}}^{\hat{p}_{2}}\circ \pi_{s+\hat{t}_{1}}^{\hat{p}_{1}}
\circ \overline{F}(\bar{t})\circ \pi_{w_{1}}^{1}(z)
$$
with the domain
$$\mathcal{D}=\Big \{(s,t,u): |s|, |t|, |u|<\overline{c}:=\frac{b}{24}
\min \big \{\hat{t}_{1}, \hat{t}_{2}, \hat{t}_{3}, c \big \} \Big\},$$
where the positive constant $b$ ($b<1$) will be determined later.
In particular, if we select $z=E^{*}_{1}$, then
$$\Psi_{E^{*}_{1}, w_{1}}(s,t, u)=\psi_{w_{1}}(s,t, u)=\psi(s,t, u),$$
which, together with (\ref{5.2.2}), implies
\begin{eqnarray}\label{5.2.4}
\mathrm{det} \Bigg(\frac{\partial \Psi_{E^{*}_{1}, w_{1}}}{\partial s}, \frac{\partial \Psi_{E^{*}_{1}, w_{1}}}{\partial t},\frac{\partial \Psi_{E^{*}_{1}, w_{1}}}{\partial u}\Bigg)\Bigg |_{(0, 0, 0)}
&\neq& 0.
\end{eqnarray}
Through a small perturbation for the function $\Psi_{z, w_{1}}(s, t, u)$,
we obtain the following function
$$\Psi_{z, w}(s,t, u)=\widetilde{F}(\tilde{t}+w_{3})\circ
\pi_{u+\hat{t}_{3}}^{\hat{p}_{3}} \circ \pi_{t+\hat{t}_{2}}^{\hat{p}_{2}}\circ \pi_{s+\hat{t}_{1}}^{\hat{p}_{1}}
\circ \overline{F}(\bar{t}+w_{2})\circ \pi_{w_{1}}^{1}(z)
$$
with the domain
$(s,t, u)\in \mathcal{D}$, $z\in \mathbb{R}^{3}_{+}$ and $w=(w_{1}, w_{2}, w')$,
where $w_{2}=(x_{1}, \ldots, x_{j})$, $w'=(y_{1}, \ldots, y_{i-4})$,
$w_{3}=(w', \widetilde{K}-\sum_{k=1}^{j}x_{k}-\sum_{k=1}^{i-4}y_{k}-w_{1}-(s+t+u))$
satisfying
$$0<w_{1}, x_{k}, y_{l}<\frac{b}{12(i+j-3)},$$
$k=1, \ldots, j$, $l=1, \ldots, i-4$, and here $\widetilde{K}=b/3$.

Note that for some fixed values of $s$, $t$, $u$ and $w$, the point
$\Psi_{z, w}(s,t, u)$ may be not necessarily accessible from
any point $z\in \mathbb{R}^{3}_{+}$ by the system (\ref{2.2}),
because two adjacent elements $\widetilde{e}_{1}$, $\widetilde{e}_{2}$ of
the ordered state set $\widetilde{\mathcal{M}}:=\{1,\bar{p}_{1},\ldots, \bar{p}_{j},
\hat{p}_{1}, \ldots, \hat{p}_{i}\}$ may be not accessible directly, i.e.,
the case where
$q_{\widetilde{e}_{1}, \widetilde{e}_{2}}=0$ may occur.
However, owing to the irreducibility of the transition rate matrix
$Q$ and the continuous dependence of the solutions on the time
and initial conditions, without loss of generality one can assume that
the ordered state set $\widetilde{\mathcal{M}}$
in the function $\Psi_{z, w}$ satisfies
that $q_{\widetilde{e}_{1}, \widetilde{e}_{2}}>0$
for any two adjacent elements $\widetilde{e}_{1}$, $\widetilde{e}_{2}\in \widetilde{\mathcal{M}}$,
i.e., the point $\Psi_{z, w}(s,t, u)$ is accessible from
any point $z\in \mathbb{R}^{3}_{+}$ by the system (\ref{2.2}).
Thus, in the remainder of the proof we will always
assume that this fact is valid for the function $\Psi_{z, w}$ considered above.

Let us define the domain
$$W=\Bigg \{x\in \mathbb{R}^{i+j-3}: 0< x_{k}< \frac{b}{12(i+j-3)}, k=1, \ldots, i+j-3 \Bigg\}.$$
By the continuous dependence of the solutions on the time
and initial conditions, it is obtained from
(\ref{5.2.4}) that there exist sufficiently small $\varepsilon>0, b$ such that
\begin{eqnarray*}
\mathrm{det} \Bigg(\frac{\partial \Psi_{z, w}}{\partial s}, \frac{\partial \Psi_{z, w}}{\partial t},\frac{\partial \Psi_{z, w}}{\partial u}\Bigg)\Bigg |_{(0, 0, 0)}
&\neq& 0
\end{eqnarray*}
for any $z\in B(E^{*}_{1}, \varepsilon)$ and $w\in W$.
Thus, for any fixed $z\in B(E^{*}_{1}, \varepsilon)$ and $w\in W$,
a small neighborhood $U_{z, w}\subset \mathcal{D}$ of the point
$(0, 0, 0)$ exists such that the mapping $\Psi_{z, w}$
is a diffeomorphism between $U_{z, w}$ and $\Psi_{z, w}(U_{z, w})$.
Note that the function $\Psi_{z, w}$ is continuously differentiable
with respect to the parameter $(z, w)$.
Due to the modified slightly inverse function theorem,
for sufficiently small $\varepsilon$ and $b$,
one can find a open subset $\overline{U}\subset\mathbb{R}^{3}_{+}$
such that for any $z\in B(E^{*}_{1}, \varepsilon)$ and $w\in W$,
we always have $U_{z, w}\subset \mathcal{D}$ and
$\overline{U}=\Psi_{z, w}(U_{z, w})$ with
$$d_{0}:= \inf_{
                  \begin{array}{c}
                    z\in B(E^{*}_{1}, \varepsilon), \\
                    w\in W, \overline{u}\in \overline{U} \\
                  \end{array}
 }
\Big |J_{z, w}(\overline{u})\Big|>0,$$
where $J_{z, w}(\overline{u})$ is the determinant of the Jacobian matrix of
$\Psi_{z, w}^{-1}$ at $\overline{u}$.

Define the event
$$\overline{\mathbf{E}}=\Bigg \{
  \begin{array}{c}
    r(0)=1, \mathrm{the}\  \mathrm{states}\  \mathrm{of}\  \mathrm{the}\  \mathrm{process}\  r(t) \
\mathrm{appear} \\
   \mathrm{ in}\  \mathrm{the}\  \mathrm{order}\  \mathrm{of}\  \{\bar{p}_{1}, \ldots, \bar{p}_{j}, \hat{p}_{1},
\ldots, \hat{p}_{i}\} \\
  \end{array}
\Bigg \}.$$
Denote by $\sigma_{0}, \sigma_{\bar{p}_{1}}, \ldots, \sigma_{\bar{p}_{j}}, \sigma_{\hat{p}_{1}}, \ldots, \sigma_{\hat{p}_{i}}$
the sojourns of the process $r(t)$ in the environmental states
$1, \bar{p}_{1}, \ldots, \bar{p}_{j}, \hat{p}_{1}, \ldots, \hat{p}_{i}$, respectively.
According to the preceding assumption, we then have
$d_{1}:= \mathbb{P}(\overline{\mathbf{E}})>0$ if $r(0)=1$.
Note that given that the event $\overline{\mathbf{E}}$
occurs, random variables $\sigma_{0}, \sigma_{\bar{p}_{1}}, \ldots, \sigma_{\bar{p}_{j}}, \sigma_{\hat{p}_{4}}, \ldots, \sigma_{\hat{p}_{i}}$
are independent. Hence,
\begin{eqnarray*}
d_{2} &:=& \mathbb{P}\Big \{\big (\sigma_{0}, \sigma_{\bar{p}_{1}}-\bar{t}_{1}, \ldots,
\sigma_{\bar{p}_{j}}-\bar{t}_{j}, \sigma_{\hat{p}_{4}}-\hat{t}_{4}, \ldots,
\sigma_{\hat{p}_{i-1}}-\hat{t}_{i-1}\big )\in W, \sigma_{\hat{p}_{i}}>\widetilde{K}
+\hat{t}_{i}\Big |\overline{\mathbf{E}} \Big \}\\[+6pt]
&=& \mathbb{P}\Big \{w\in W, \sigma_{\hat{p}_{i}}>\widetilde{K}
+\hat{t}_{i}\Big |\overline{\mathbf{E}} \Big \}>0.
\end{eqnarray*}
Let $K=\sum_{k=1}^{j}\bar{t}_{k}+\sum_{k=1}^{i}\hat{t}_{k}+\widetilde{K}$.
For any $z\in B(E^{*}_{1}, \varepsilon)$ and any Borel set $B\subset \overline{U}$, we have
\begin{eqnarray*}
&&\mathbb{P}\big(K, z, 1, B \times \{\hat{p}_{i}\}\big) \\[+6pt]
&\geq& d_{1}\cdot \mathbb{P}\Bigg \{
\begin{array}{c}
\big (\sigma_{0}, \sigma_{\bar{p}_{1}}-\bar{t}_{1}, \ldots,
\sigma_{\bar{p}_{j}}-\bar{t}_{j}, \sigma_{\hat{p}_{4}}-\hat{t}_{4}, \ldots,
\sigma_{\hat{p}_{i-1}}-\hat{t}_{i-1}\big )=w \in W, \\
                                  \big (\sigma_{\hat{p}_{1}}-\hat{t}_{1}, \sigma_{\hat{p}_{2}}-\hat{t}_{2}, \sigma_{\hat{p}_{3}}-\hat{t}_{3} \big )=(s,t,u)\in \Psi_{z, w}^{-1}(B), \sigma_{\hat{p}_{i}}>\widetilde{K}
+\hat{t}_{i}\\
                                \end{array}
\Bigg |\overline{\mathbf{E}} \Bigg\} \\[+6pt]
&=& d_{1}\cdot \mathbb{P} \Big \{ (s,t,u)\in \Psi_{z, w}^{-1}(B)\Big | \overline{\mathbf{E}},
w\in W, \sigma_{\hat{p}_{i}}>\widetilde{K}
+\hat{t}_{i}\Big \}\cdot \mathbb{P}\Big\{ w\in W, \sigma_{\hat{p}_{i}}>\widetilde{K}
+\hat{t}_{i} \Big |\overline{\mathbf{E}} \Big\}\\[+6pt]
&\geq& d_{1}d_{2}\cdot  \mathbb{P}\Big \{ (s,t,u)\in \Psi_{z, w}^{-1}(B)\Big | \overline{\mathbf{E}},
w\in W, \sigma_{\hat{p}_{i}}>\widetilde{K}
+\hat{t}_{i}\Big \}.
\end{eqnarray*}
Let $h(x_{1},x_{2},x_{3})$ be the probability density function of random variables
$(\sigma_{\hat{p}_{1}}, \sigma_{\hat{p}_{2}}, \sigma_{\hat{p}_{3}})$
under the condition that the event $\overline{\mathbf{E}}$
occurs. Since $\sigma_{\hat{p}_{1}}$, $\sigma_{\hat{p}_{2}}$, $\sigma_{\hat{p}_{3}}$
are independent exponential random variables given that
the event $\overline{\mathbf{E}}$
occurs, the function $h(s+\hat{t}_{1}, t+\hat{t}_{2}, u+\hat{t}_{3})$
is a smooth function and
$d_{3}:= \inf_{(s,t,u)\in \mathcal{D}} h(s+\hat{t}_{1}, t+\hat{t}_{2}, u+\hat{t}_{3})>0$.
Thus,
\begin{eqnarray}\label{5.2.5}
\mathbb{P}\big(K, z, 1, B \times \{\hat{p}_{i}\}\big)
&\geq& d_{1}d_{2}\cdot \int_{U_{z,w}}h(s+\hat{t}_{1}, t+\hat{t}_{2}, u+\hat{t}_{3})\mathbf{1}_{\Psi_{z, w}^{-1}(B)}(s,t,u)dsdtdu \nonumber\\[+6pt]
&\geq& d_{1}d_{2}d_{3}\cdot \int_{U_{z,w}}\mathbf{1}_{\Psi_{z, w}^{-1}(B)}(s,t,u)dsdtdu\nonumber\\[+6pt]
&=& d_{1}d_{2}d_{3}\cdot \int_{\overline{U}}\mathbf{1}_{B}(u_{1}, u_{2} ,u_{3})\Big |J_{z, w}(u_{1}, u_{2} ,u_{3})\Big|du_{1}du_{2}du_{3} \nonumber\\[+6pt]
&\geq& d_{0}d_{1}d_{2}d_{3}\cdot \widetilde{m}(B)=d_{4}\cdot m( B \times \{\hat{p}_{i}\}),
\end{eqnarray}
where $d_{4}$ is some positive constant.
By Lemma \ref{lemma4.3.} and the similar arguments in Step 2 of the proof
in Theorem \ref{theorem4.1.}, it is obtained that for any initial
value $(z_{0}, e)=((S(0), I(0), R(0)), e)\in X$,
there exists with probability one a sequence of strictly increasing
finite stopping times $\{\varsigma_{n}\}_{n=1}^{\infty}$
with respect to filtration $\mathcal{F}_{0}^{n}$,
such that $r(\tau_{\varsigma_{n}})=1$ and $(S(\tau_{\varsigma_{n}}), I(\tau_{\varsigma_{n}}), R(\tau_{\varsigma_{n}}))\in
B(E^{*}_{1}, \varepsilon)$ for all $n\in \mathbb{N}$.
Combining this property and (\ref{5.2.5}), one can obtain by the strong Markov property
of the process $((S(t), I(t), R(t)), r(t))$ that
for $d:=\mathbb{P}\big(K, z, 1, B \times \{\hat{p}_{i}\}\big)>0$,
\begin{eqnarray*}
\mathbb{P}(\widetilde{\mathbf{E}})
&\geq& d+(1-d)d+(1-d)^{2}d+\cdots \\[6pt]
&=&\frac{ d}{1-(1-d)}=1,
\end{eqnarray*}
where $\widetilde{\mathbf{E}}$ represents
the event that the process $((S(t), I(t), R(t)), r(t))$
will enter $B \times \{\hat{p}_{i}\}$ at some finite moment.
Hence, if $m( B \times \{\hat{p}_{i}\})>0$
with $B\subset \overline{U}$, then $\mathbb{P}(\widetilde{\mathbf{E}})=1$,
which implies that the condition (a) in the definition of Harris recurrent
is satisfied with the measure $\phi_{1}(A)=m\big(A\cap (B \times \{\hat{p}_{i}\})\big)$
for any $A\in \mathcal{B}(X)$. Note that from Theorem \ref{theorem5.1.},
we have got that the process $((S(t), I(t), R(t)), r(t))$
has an invariant probability measure $\nu^{*}$
on the state space $X$. Consequently, the process $((S(t), I(t), R(t)), r(t))$
is positive Harris recurrent.

{\bf Step 2.} Next, we shall show that
a lattice distribution $\mathfrak{a}$ exists such that
the corresponding kernel $\mathbf{K}_{\mathfrak{a}}$ admits an
everywhere non-trivial continuous component,
which implies that
the stochastic process $((S(t), I(t), R(t)), r(t))$
is a $\mathcal{T}$-process. Combining with the result in Step 1, it then follows from
Proposition \ref{corollary5.1.} that for any initial value $x_{0}\in X$,
the assertions (\ref{theorem19-1}) and (\ref{theorem19-2}) are valid.
Let us first consider two cases as follows.

Case 1. $r(0)=1\in \mathcal{M}$. For each $k\in \mathbb{N}$,
it is obtained from Lemma \ref{lemma-global} and the globally asymptotic
stability of the equilibrium $E^{*}_{1}$ of system (\ref{2.2}) in the state 1
that there exists some $n_{k}^{1}\in \mathbb{N}$ such that
if $n\geq n_{k}^{1}$, $\pi_{n K}^{1}(z)\in B(E^{*}_{1}, \varepsilon)$
for all $z\in [k^{-1}, k]^{3}\cap \mathcal{K}$.

Case 2. $r(0)=e \in \mathcal{M}$ but $e\neq 1$.
Due to the connectivity between distinct states of the Markov
chain $r(t)$, we can get that a positive integer $l_{e}$ ($1\leq l_{e}\leq E-1$)
exists such that
$$q_{h_{0}^{e}, h_{1}^{e}}q_{h_{1}^{e}, h_{2}^{e}}\cdots q_{h_{l_{e}-1}^{e}, h_{l_{e}}^{e}} >0,$$
where $\{h_{k}^{e}: 0\leq k\leq l_{e}\}\subset \mathcal{M}$ and $h_{0}^{e}=e$, $h_{l_{e}}^{e}=1$.
From Lemma \ref{lemma-global} and the globally asymptotic
stability of the equilibrium $E^{*}_{1}$,
it then follows that for each $k\in \mathbb{N}$,
some $n_{k}^{e}\in \mathbb{N}$ exists such that if $n\geq n_{k}^{e}$,
$$\pi_{nK-(s_{0}+\cdots +s_{l_{e}-1})}^{h_{l_{e}}^{e}}\circ \pi_{s_{l_{e}-1}}^{h_{l_{e}-1}^{e}}
\circ \cdots \pi_{s_{0}}^{h_{0}^{e}}(z)\in B(E^{*}_{1}, \varepsilon)$$
for all $z\in [k^{-1}, k]^{3}\cap \mathcal{K}$ and
$(s_{0}, \cdots, s_{l_{e}-1})\in [0, K/l_{e}]^{l_{e}}$.

Let us select a sequence of strictly increasing $\{n_{k}\}_{k=1}^{\infty}$
with $n_{k}=\max_{e\in \mathcal{M}}\{n_{k}^{e}\}$.
Denote by $\sigma_{h_{0}^{e}}, \ldots, \sigma_{h_{l_{e}}^{e}}$ the sojourn times of the process $r(t)$
in the states $h_{0}^{e}, \ldots, h_{l_{e}}^{e}$, respectively.
For each $k\in \mathbb{N}$, we then have:

(1) if $r(0)=1$,
$$p_{k}^{1}=\mathbb{P}\Big\{\tilde{\sigma}>n_{k}K \Big| r(0)=1\Big\}>0,$$
where $\tilde{\sigma}$ denotes the sojourn time of the process $r(t)$
in the states 1;

(2) if $r(0)=e\in \mathcal{M}\backslash\{1\}$,
$$p_{k}^{e}=\mathbb{P}\Big \{ (\sigma_{h_{0}^{e}}, \ldots, \sigma_{h_{l_{e}-1}^{e}})\in [0, K/l_{e}]^{l_{e}}, \sigma_{h_{0}^{e}}+\cdots+\sigma_{h_{l_{e}}^{e}}> n_{k}K \Big| r(0)=e \Big\}>0.$$

Letting $p_{k}=\min_{e\in \mathcal{M}}\{p_{k}^{e}\}$ for each $k\in \mathbb{N}$, we can get
$$\mathbb{P}\Big (n_{k}K, z, e, B(E^{*}_{1}, \varepsilon)\times \{1\}\Big)
\geq p_{k}$$
for any $z\in [k^{-1}, k]^{3}\cap \mathcal{K}$ and $e\in \mathcal{M}$.
By the Kolmogorov-Chapman equation, this, together with (\ref{5.2.5}), implies that
$$\mathbb{P}\Big ((n_{k}+1)K, z, e, B\times \{1\}\Big)
\geq p_{k}d_{4}\cdot m(B\times \{1\})$$
for any $z\in [k^{-1}, k]^{3}\cap \mathcal{K}$, $e\in \mathcal{M}$
and any Borel set $B\subset \overline{U}$,
where we here use the fact that in the proof of Step 1 in this theorem,
we can assume that the state $\hat{p}_{i}=1$ without loss of generality.
As in the proof of Theorem 4.2 in \cite{Yingang2014},
we construct a lattice distribution $\mathfrak{a}(nK)=2^{-n}$, $n\in \mathbb{N}$,
and the corresponding Markov transition function
$\mathbf{K}_{\mathfrak{a}}(z, e, A)=\sum_{n=1}^{\infty}2^{-n}\mathbb{P}(nK, z,e,A)$
for any $(z, e)\in X$ and $A\in \mathcal{B}(X)$.
Moreover, the kernel $\mathbf{K}_{\mathfrak{a}}$ has an
everywhere non-trivial continuous component
$\mathcal{T}: X\times \mathcal{B}(X)\rightarrow [0, +\infty)$
defined by
$$\mathcal{T}(z,e,A)=2^{-(n_{k+1}+1)}p_{k+1}d_{4}\cdot
m\Big(A \cap (\overline{U} \times \{1\})\Big)$$
when $z\in \big(((k+1)^{-1}, k+1)^{3}\backslash (k^{-1}, k)^{3} \big)\cap \mathcal{K}$,
$k\in \mathbb{N}$. Thus,
the process $((S(t), I(t), R(t)), r(t))$
is a $\mathcal{T}$-process.

{\bf Step 3.} We now shall that the stationary distribution $\nu^{*}$
of the process $((S(t)$, $I(t)$, $R(t))$, $r(t))$ has a density $f^{*}$
with respect to the product measure $m$ on $X$
and $\mathrm{supp}(f^{*})= \Gamma\times \mathcal{M}$.
Since the argument is similar to that of \cite{Yingang2011},
we here only sketch the proof to point out the difference with it.

Noting that $\Gamma$ is a positive invariant set of system (\ref{2.2}),
it follows from (b) of Theorem \ref{theorem4.1.} that
$\lim_{t\rightarrow +\infty}\mathbb{P}\{z(t)\in \Gamma\}=1$,
which implies that $\nu^{*}(\Gamma\times \mathcal{M})=1$.
Applying the Lebesgue decomposition theorem, there exist unique
probability measures $\nu^{*}_{a}$, $\nu^{*}_{s}$ and $\kappa\in[0, 1]$
such that $\nu^{*}=(1-\kappa)\nu^{*}_{a}+ \kappa \nu^{*}_{s}$, where $\nu^{*}_{a}$
is absolutely continuous with respect to $m$ and $\nu^{*}_{s}$
is singular. It is easy to see that if $\kappa=0$, then $\nu^{*}=\nu^{*}_{a}$,
i.e., $\nu^{*}$ is absolutely continuous with respect to $m$.
Suppose $\kappa\neq0$. By the similar arguments in Proposition 3.1 of \cite{Yingang2011},
one can get that $\nu^{*}_{s}$ is also a stationary distribution,
and a measurable subset $\Gamma_{0}\subset \Gamma$ exists
with $\widetilde{m}(\Gamma_{0})=0$ and $\nu^{*}_{s}(\Gamma_{0}\times \mathcal{M})=1$.
From the proof in (b) of Theorem \ref{theorem4.1.} and Step 1 of the proof in this theorem,
it is obtained by the continuous dependence of the solutions
on the time and initial conditions that for sufficiently
small $\varepsilon$, $b$ in the Step 1, a open set $\overline{U}\subset \Gamma$
exists such that
$$\mathbb{P}\Big(K, z, 1, B\times \{\hat{p}_{i}\}\Big)\geq d_{4}\cdot m(B\times \{\hat{p}_{i}\})$$
for any Borel set $B \subset \overline{U}$,
where $z\in B(E^{*}_{1}, \varepsilon)$.
Noting that $m\big((A\backslash C)\times \overline{\mathcal{M}}\big)= m(A\times \overline{\mathcal{M}})$
for any $\overline{\mathcal{M}} \subset \mathcal{M}$
and any Lebesgue measurable subsets $A$, $C$ of $\mathcal{K}$
with $\widetilde{m}(C)=0$, we then have
\begin{eqnarray*}
\mathbb{P}\Big( K, z, 1, (\Gamma \backslash \Gamma_{0})\times \mathcal{M}\Big)
&\geq& \mathbb{P}\Big( K, z, 1, (\overline{U} \backslash \Gamma_{0})\times \mathcal{M}\Big)  \\[+4pt]
&\geq& \mathbb{P}\Big( K, z, 1, (\overline{U} \backslash \Gamma_{0})\times \{\hat{p}_{i}\}\Big) \\[+4pt]
&\geq& d_{4}\cdot m\big( (\overline{U} \backslash \Gamma_{0})\times \{\hat{p}_{i}\}\big) \\[+4pt]
&=& d_{4}\cdot m \big(\overline{U} \times \{\hat{p}_{i}\}\big)>0.
\end{eqnarray*}
Through the similar arguments in Proposition 3.1 in \cite{Yingang2011},
it follows that $\nu^{*}_{s}\mathbb{P}^{K}\big((\Gamma \backslash \Gamma_{0})\times \mathcal{M}\big)>0$,
which is contradict with
$$\nu^{*}_{s}\mathbb{P}^{K}\big((\Gamma \backslash \Gamma_{0})\times \mathcal{M}\big)
=\nu^{*}_{s}\big((\Gamma \backslash \Gamma_{0})\times \mathcal{M}\big)=0$$
due to the stationary distribution property. Thus, $\kappa=0$, which implies that
$\nu^{*}$ is absolutely continuous with respect to $m$
with the density $f^{*}$.

Finally, by the similar arguments of Proposition 3.1 in \cite{Yingang2011},
one can also conclude that
$\int_{V}f^{*}dm>0$ for all open set $V\subset \Gamma \times \mathcal{M}$, which implies that
$\mathrm{supp}(f^{*})= \Gamma\times \mathcal{M}$.
This completes the proof of Theorem \ref{theorem5.2.}.\ \ $\Box$

\end{pf}

\section{Discussion}
\label{Conclusion}

In this paper, we have established the threshold dynamics
of the disease extinction and persistence
for the system (\ref{2.2}).
That is,
the disease can be eradicated almost surely if $\mathcal{R}_{0}<1$,
while the disease persists almost surely if $\mathcal{R}_{0}>1$.
Moreover,
we also have given the ergodic analyzing of this epidemic model.
The global attractivity of the $\Omega$-limit set of the system
and the convergence in total variation of the
instantaneous measure to the stationary measure
were obtained under the weakened conditions.
From two environmental regimes to any finite ones,
the method presented in this paper is a generalization of
the techniques used by \cite{Yingang2011,Yingang2014} for
the ergodic analyzing of the piecewise deterministic Markov process.
For analyzing the ergodicity of the considered system,
the generalized method requires weaker conditions and
is applicable to the multi-dimension system with
any number of environmental regimes.

Note also that only the transmission rate $\beta$ of model (\ref{2.2}) is
disturbed because in reality it is more sensitive to environmental
fluctuations than other parameters
for human populations. Nevertheless,
it is easy to see that the method used in this paper can
be easily extended to the case where
other parameters
of model (\ref{2.2}), such as the recruitment rate $\Lambda$
and the natural mortality $\mu$, can also change with
the switching of environmental regimes,
which may be more reasonable for the wildlife population.
In addition, by this new method, under weaker conditions,
we can directly extend the results
in \cite{Yingang2011,Yingang2014,HIEU2015}
from two environmental states to any finite ones.

However, in the case when $\mathcal{R}_{0}>1$,
the extra condition that some point in the $\Omega$-limit subset $\Gamma$
satisfying
the condition (\textbf{H}) exists is required for ensuring
the global attractivity of the $\Omega$-limit set of system (\ref{2.2})
and the convergence in total variation of the
instantaneous measure to the stationary measure.
From a biological point of view, what is the biological meaning of this condition?
To obtain the ergodicity of system (\ref{2.2}),
is it enough that only the condition that $\mathcal{R}_{0}>1$ is satisfied?
We leave these questions for future work.

\section*{Acknowledgments}  The authors would extend their
thanks to Professors Jianhua Huang
and Weiming Wang for their valuable comments.
D. Li and S. Liu are supported by  the National Natural
Science Foundation of China (No.11471089).
J. Cui is supported by the National Natural
Science Foundation of China (No.11371048).

\end{document}